\documentclass[mnsc,nonblindrev]{informs3_nojournal}

\OneAndAHalfSpacedXII 

\usepackage[utf8]{inputenc}
\usepackage{algorithmic}
\usepackage{algorithm}
\usepackage{subcaption}
\usepackage{csquotes}
\usepackage[section]{placeins}

\usepackage{enumitem}

\usepackage[T1]{fontenc}    
\usepackage{url}            
\usepackage{booktabs}       
\usepackage{amsfonts}       
\usepackage{nicefrac}       

\newcommand{\E} { \mathbb{E} }

\newcommand{\trace}{\text{trace}}

\newcommand{\bbR}{\mathbb{R}}

\newcommand{\bbE}{\mathbb{E}}
\newcommand{\bbP}{\mathbb{P}}

\newcommand{\calX}{\mathcal{X}}
\newcommand{\calH}{\mathcal{H}}

\newcommand{\calD}{\mathcal{D}}
\newcommand{\supp}{\xi_S}
\newcommand{\ellTw}{\ell^{w^\ast}_{\mathrm{SPO}}}
\newcommand{\ellT}{\ell_{\mathrm{SPO}}}
\newcommand{\riskT}{R_{\mathrm{SPO}}}
\newcommand{\ellS}{\ell_{\mathrm{SPO+}}}

\newcommand{\riskS}{R_{\mathrm{SPO+}}}
\newcommand{\LSPO}{L^n_{\mathrm{SPO+}}}
\newcommand{\degg}{\text{deg}}

\newcommand{\edit}{\textcolor{black}}
\newcommand{\blockedit}{\color{black}}

\usepackage{natbib}
 \bibpunct[, ]{(}{)}{,}{a}{}{,}%
 \def\bibfont{\small}%
 \def\bibsep{\smallskipamount}%
 \def\bibhang{24pt}%
\TheoremsNumberedThrough     
\ECRepeatTheorems

\EquationsNumberedThrough    

\MANUSCRIPTNO{}

\begin{document}


\RUNAUTHOR{Elmachtoub and Grigas}

\RUNTITLE{Smart \enquote{Predict, then Optimize}}

\TITLE{Smart \enquote{Predict, then Optimize}}

\ARTICLEAUTHORS{%
\AUTHOR{Adam N. Elmachtoub}
\AFF{Department of Industrial Engineering and Operations Research and Data Science Institute, Columbia University, New York, NY 10027, \EMAIL{adam@ieor.columbia.edu}}
\AUTHOR{Paul Grigas}
\AFF{Department of Industrial Engineering and Operations Research, University of California, Berkeley, CA 94720, \EMAIL{pgrigas@berkeley.edu}}
} 

\ABSTRACT{Many real-world analytics problems involve two significant challenges:  prediction and optimization. Due to the typically complex nature of each challenge, the standard paradigm is predict-then-optimize. By and large, machine learning tools are intended to minimize prediction error and do not account for how the predictions will be used in the downstream optimization problem. In contrast, we propose a new and very general framework, called Smart ``Predict, then Optimize'' (SPO), which directly leverages the optimization problem structure, i.e., its objective and constraints, for designing better prediction models.  A key component of our framework is the SPO loss function which measures the decision error induced by a prediction. 

Training a prediction model with respect to the SPO loss is computationally challenging, and thus we derive, using duality theory, a convex surrogate loss function which we call the SPO+ loss. Most importantly, we prove that the SPO+ loss is statistically consistent with respect to the SPO loss under mild conditions. Our SPO+ loss function can tractably handle any polyhedral, convex, or even mixed-integer optimization problem with a linear objective. Numerical experiments on shortest path and portfolio optimization problems show that the SPO framework can lead to significant improvement under the predict-then-optimize paradigm, in particular when the prediction model being trained is misspecified. 
We find that linear models trained using SPO+ loss tend to dominate random forest algorithms, even when the ground truth is highly nonlinear. 
}%



\maketitle
\vspace{-25pt}

\section{Introduction}\label{sec:intro}

In many real-world analytics applications of operations research, a combination of both machine learning and optimization are used to make decisions. Typically, the optimization model is used to generate decisions, while a machine learning tool is used to generate a prediction model that predicts key unknown parameters of the optimization model. Due to the inherent complexity of both tasks,  a broad purpose approach that is often employed in analytics practice is the \textit{predict-then-optimize} paradigm. 

For example, consider a vehicle routing problem that may be solved several times a day.  First, a previously trained prediction model provides predictions for the travel time on all edges of a road network based on current traffic, weather, holidays, time, etc. Then, an optimization solver provides near-optimal routes using the predicted travel times as input.  We emphasize that most solution systems for real-world analytics problems involve some component of both prediction and optimization  (see \cite{angalakudati2014business,chan2012optimizing,deo2015planning,gallien2015initial,cohen2017impact,besbes2015optimization,mehrotra2011or,chan2013prioritizing,ferreira2015analytics} for recent examples and 
recent expositions by \cite{simchi2013om, denbridging,deng2018coalescing,mivsic2020data}). 
Except for a few limited options, machine learning tools do not effectively account for how the predictions will be used in a downstream optimization problem. In this paper, we provide a general framework called  Smart ``Predict, then Optimize'' (SPO) for training prediction models that effectively utilize the structure of the nominal optimization problem, i.e., its constraints and objective. \textit{Our SPO framework is fundamentally designed to generate prediction models that aim to minimize decision error, not prediction error.} 

One key benefit of our SPO approach is that it maintains the decision paradigm of sequentially predicting and then optimizing.  However, when training our prediction model, the structure of the nominal optimization problem is explicitly used. 
The quality of a prediction is \textit{not} measured based on prediction error such as least squares loss or other popular loss functions. Instead, in the SPO framework, the quality of a prediction is measured by the decision error. That is, suppose a prediction model is trained using historical feature data $(x_1,\ldots, x_n)$ and associated parameter data $(c_1,\ldots,c_n)$. Let $(\hat{c}_1,\ldots, \hat{c}_n)$ denote the predictions of the parameters under the trained model. The least squares (LS) loss, for example, measures error with the squared norm $\|c_i - \hat{c}_i \|_2^2$, completely ignoring the decisions induced by the predictions. In contrast, the \textit{SPO loss} is the true cost of the decision induced by $\hat{c}_i$ minus the optimal cost under the true parameter $c_i$. In the context of vehicle routing, the SPO loss measures the extra travel time incurred due to solving the routing problem on the predicted, rather than true, edge cost parameters.


In this paper, we focus on predicting unknown parameters of a contextual stochastic optimization problem, where the parameters appear linearly in the objective function, i.e., the cost vector of any linear, convex, or integer optimization problem. The core of our SPO framework is a new loss function for training prediction models.
Since the SPO loss function is difficult to work with, significant effort revolves around deriving a surrogate loss function, SPO+, that is convex and therefore can be optimized efficiently. To show the validity of the surrogate SPO+ loss, we prove a highly desirable statistical consistency property, and show it performs well empirically  compared to standard predict-then-optimize approaches. {\blockedit In essence, we prove that the function that minimizes the Bayes risk associated to the SPO+ loss is the regression function $\mathbb{E}[c | x]$, which also minimizes the Bayes risk of the SPO loss (under mild assumptions). Interestingly, $\mathbb{E}[c | x]$ also minimizes the Bayes risk associated to the LS loss under the same conditions. Thus, SPO+ and LS (or any convex combination of the two) are essentially on ``equal footing'' -- they are both theoretically valid (consistent) and computationally tractable choices for the loss function. However, when the ultimate goal is to solve a downstream optimization task, the SPO+ loss is the natural choice as it is tailored to the optimization problem and works significantly better in practice than LS.  

Empirically, we observe that even when the prediction task is challenging due to model misspecification, the SPO framework can still yield near-optimal decisions. We note that a fundamental property of the SPO framework is the requirement that the prediction is directly ``plugged in'' to the downstream optimization problem. An alternative procedure may alter the decision making process in some way, such as by adding robustness or by taking into account the entire dataset (instead of just the prediction). A strong advantage of our SPO approach is that it has good performance even when the naive prediction problem is challenging, see the illustrative example in Section \ref{sect:illustrative}. Another advantage is that the downstream optimization problem is typically more computationally tractable and more attractive to practitioners than a more complex alternative procedure. 
On the other hand, alternative decision making procedures may provide other advantages, such as improved generalization performance via the introduction of bias and/or robustness. However, designing such procedures is  more challenging in the presence of contextual data and  combining them with the SPO approach would be worthwhile of future research.
Overall, we believe our SPO framework provides a clear foundation for designing operations-driven machine learning tools that can  be leveraged in real-world optimization settings. }

\setlist{nolistsep}
Our contributions may be summarized as follows:
\begin{enumerate}[leftmargin=*]
\item We first formally define a new loss function, which we call the SPO loss, that measures the error in predicting the cost vector of a nominal optimization problem with linear, convex, or integer constraints. 
The loss corresponds to the suboptimality gap -- with respect to the true/historical cost vector -- due to implementing a possibly incorrect decision induced by the predicted cost vector.
Unfortunately, the SPO loss function can be nonconvex and discontinuous in the predictions, implying that training ML models under the SPO loss may be challenging. 

\item Given the intractability of the SPO loss function, we develop a surrogate loss function which we call the SPO+ loss. This surrogate loss function is derived using a sequence of steps motivated by duality theory (Proposition \ref{alternative-rep}), a data scaling approximation, and a first-order approximation.
The resulting SPO+ loss function is convex in the predictions (Proposition~\ref{prop:convex_ub}), which allows us to design an algorithm based on stochastic gradient descent for minimizing SPO+ loss (Proposition~\ref{prop:phi_subgrad}). Moreover, when training a linear regression model to predict the objective coefficients of a linear program, only a linear optimization problem needs be solved to minimize the SPO+ loss (Proposition~\ref{reformulation}). 

\item We prove a fundamental connection to classical machine learning under a very simple and special instance of our SPO framework. Namely, under this instance the SPO loss is exactly the 0-1 classification loss (Proposition \ref{binary1}) and the SPO+ loss is exactly the hinge loss (Proposition \ref{binary2}). The hinge loss is the basis of the popular SVM method and is a surrogate loss to approximately minimize the 0-1 loss, and thus our framework generalizes this concept to a very wide family of optimization problems with constraints. 


\item We prove a key consistency result of the SPO+ loss function (Theorem \ref{thm:main}, Proposition \ref{true_consistency}, Proposition \ref{spo_consistent}), which further motivates its use. Namely, under full distributional knowledge, minimizing the SPO+ loss function is in fact equivalent to minimizing the SPO loss if two mild conditions hold:  the distribution of the cost vector (given the features) is continuous and symmetric about its mean. For example, these assumptions are satisfied by the standard Gaussian noise approximation. \textit{This consistency property is widely regarded as an essential property of any surrogate loss function across the statistics and machine learning literature.} For example, the famous hinge loss and logistic loss functions are consistent with the 0-1 classification loss. 

\item Finally, we validate our framework through numerical experiments on the shortest path and portfolio optimization problem. We test our SPO framework against standard predict-then-optimize approaches, and evaluate the out of sample performance with respect to the SPO loss. Generally, the value of our SPO framework increases as the degree of model misspecification increases. This is precisely due to the fact the SPO framework makes ``better'' wrong predictions, essentially ``tricking'' the optimization problem into finding near-optimal solutions. 
Remarkably, a linear model trained using SPO+ even dominates a state-of-the-art random forests algorithm, even when the ground truth is highly nonlinear. 
\end{enumerate}

\subsection{Applications} \label{apps}
Settings where the input parameters (cost vectors) of an optimization problem need to be predicted from contextual (feature) data are numerous. Let us now highlight a few, of potentially many, application areas for the SPO framework.

\paragraph{Vehicle Routing.} In numerous applications, the cost of each edge of a graph needs to be predicted before making a routing decision. The cost of an edge typically corresponds to the expected length of time a vehicle would need to traverse the corresponding edge. For clarity, let us focus on one important example -- the shortest path problem.  In the shortest path problem, one is given a weighted directed graph, along with an origin node and destination node, and the goal is to find a sequence of edges from the origin to the destination at minimum possible cost. A well-known fact is that the shortest path problem can be formulated as a linear optimization problem, but there are also alternative specialized algorithms such as the famous Dijkstra's algorithm (see, e.g., \cite{ahuja1993network}). The data used to predict the cost of the edges may incorporate the length, speed limit, weather, season, day, and real-time data from mobile applications such as Google Maps and Waze. Simply minimizing prediction error may not suffice nor be appropriate, as over- or under-predictions have starkly different effects across the network. The SPO framework would ensure that the predicted weights lead to shortest paths, and would naturally emphasize the estimation of edges that are critical to this decision. See Figure \ref{fig:example} in Section \ref{sec:framework} for an in-depth example. 

\paragraph{Inventory Management.} 
In inventory planning problems such as the economic lot sizing problem (\cite{wagner1958dynamic}) or the joint replenishment problem (\cite{levi2006primal}), the demand is the key input into the optimization model. In practical settings, demand is highly nonstationary and can depend on historical and contextual data such as weather, seasonality, and competitor sales. The decisions of when to order inventory are captured by a linear or integer optimization model, depending on the complexity of the problem. Under a common formulation (see \cite{levi2006primal,cheung2016submodular}), the demand appears linearly in the objective, which is convenient for the SPO framework. The goal is to design a prediction model that maps feature data to demand predictions, which in turn lead to good inventory plans. 

\paragraph{Portfolio Optimization.} 
In financial services applications, the returns of potential investments need to be somehow estimated from data, and can depend on many features which typically include historical returns, news, economic factors, social media, and others. In portfolio optimization, the goal is to find a portfolio with the highest return subject to a constraint on the total risk, or variance, of the portfolio. While the returns are often highly dependent on auxiliary feature information, the variances are typically much more stable and are not as difficult nor sensitive to predict. Our SPO framework would result in predictions that lead to high performance investments that satisfy the desired level of risk. A least squares loss approach places higher emphasis on estimating higher valued investments, even if the corresponding risk may not be ideal. In contrast, the SPO framework directly accounts for the risk of each investment when training the prediction model. 

\subsection{Related Literature} \label{literature}
Perhaps the most related work is that of \citet{kao2009directed}, who also directly seek to train a machine learning model that minimizes loss with respect to a nominal optimization problem. In their framework, the nominal problem is an unconstrained quadratic optimization problem, where the unknown parameters appear in the linear portion of the objective. 
Their work does not extend to settings where the nominal optimization problem has constraints, which our framework does. \cite{donti2017task} proposes a heuristic to address a more general setting than that of \citet{kao2009directed}, and also focus on the case of quadratic optimization. 
These works also bypass issues of non-uniqueness of solutions of the nominal problem (since their problem is strongly convex), which must be addressed in our setting to avoid degenerate prediction models. 

In \citet{ban2019big}, ML models are trained to directly predict the optimal solution of a newsvendor problem from data. Tractability and statistical properties of the method are shown as well as its effectiveness in practice. However, it is not clear how this approach can be used when there are constraints, since feasibility issues may arise. 

The general approach in \citet{bertsimas2020predictive} considers the problem of accurately estimating an unknown optimization objective using machine learning models, specifically ML models where the predictions can be described as a weighted combination of training samples, e.g., nearest neighbors and decision trees. In their approach, they estimate the objective of an instance by applying the same weights generated by the ML model to the corresponding objective functions of those samples. This approach differs from standard predict-then-optimize \textit{only} when the objective function is nonlinear in the unknown parameter. Note that the unknown parameters of all the applications mentioned in Section \ref{apps} appear linearly in the objective. Moreover, the training of the ML models does not rely on the structure of the nominal optimization problem, in contrast to the SPO framework.

The approach in \citet{tulabandhula2013machine} relies on minimizing a loss function that combines the prediction error with the operational cost of the model on an unlabeled dataset. However, the operational cost is with respect to the predicted parameters, and not the true parameters. \cite{gupta2017small} consider combining estimation and optimization in a setting without features/contexts. We also note that our SPO loss, while mathematically different, is similar in spirit to the notion of relative regret introduced in \cite{lim2012robust} in the specific context of portfolio optimization with historical return data and without features.    
Other approaches for finding near-optimal solutions from data include operational statistics (\cite{liyanage2005practical,chu2008solving}), sample average approximation (\cite{kleywegt2002sample,schutz2009supply,bertsimas2018robust}), and robust optimization (\cite{bertsimas2006robust,bertsimas2018data,wang2016likelihood}). There has also been some recent progress on submodular optimization from samples (\cite{balkanski2016power,balkanski2017limitations}). These approaches typically do not have a clear way of using feature data, nor do they directly consider how to train a machine learning model to predict optimization parameters.  

Another related stream of work is in data-driven inverse optimization, where feasible or optimal solutions to an optimization problem are observed and the objective function has to be learned (\cite{aswani2018inverse,keshavarz2011imputing, chan2014generalized, bertsimas2015data,esfahani2018data}).  In these problems, there is typically a single unknown objective, and no previous samples of the objective are provided. We also note there have been recent approaches for regularization (\cite{ban2018machine}) and model selection (\cite{besbes2010testing,den2016decision,sen2017learning}) in the context of an optimization problem.  


Lastly, we note that our framework is related to the general setting of structured prediction (see, e.g., \citet{taskar2005learning, tsochantaridis2005large, nowozin2011structured, osokin2017structured} and the references therein). Motivated by problems in computer vision and natural language processing, structured prediction is a version of multiclass classification that is concerned with predicting structured objects, such as sequences or graphs, from feature data.
The SPO+ loss is similar in spirit to that of the structured SVM (SSVM) and is indeed a convex, upper bound on the SPO loss, akin to the SSVM. However, there are fundamental differences with our approach and the the SSVM approach. In the SSVM approach, the structured object one would be predicting is the decision $w$ directly from the feature $x$ (\cite{taskar2005learning}). In our setting, we have access to historical data on $c$ which is richer than observations of decisions, since cost vectors induce optimal decisions naturally. Under one special case of our framework, we prove that the SPO loss is equivalent to 0/1 loss, while the SPO+ loss is equivalent to the hinge loss. Thus, our framework can be seen as a type of generalization of the SSVM. Finally, we remark that our derivation of the surrogate SPO+ loss relies on completely new ideas using duality theory, which help explain the strong empirical performance.

\section{``Predict, then Optimize'' Framework} 
\label{sec:framework}
We now describe the ``Predict, then Optimize'' framework which is central to many applications of optimization in practice. Specifically, we assume that there is a nominal optimization problem of interest with a linear objective, where the decision variable $w \in \bbR^d$ and feasible region  $S \subseteq \bbR^d$ are well-defined and known with certainty. However, the cost vector of the objective, $c \in \bbR^d$, is not available at the time the decision must be made; instead, an associated feature vector $x \in \bbR^p$  is available. Let $\mathcal{D}_x$ be the conditional distribution of $c$ given $x$. The goal for the decision maker is to solve, for any new instance characterized by $x$, is to solve the contextual stochastic optimization problem 
\begin{align}
\min_{w \in S} \E_{c \sim \mathcal{D}_x}[c^\top w |x] \ = \ \min_{w \in S} \E_{c \sim \mathcal{D}_x}[c|x]^\top w \ .
\label{core}
\end{align}

The predict-then-optimize framework relies on using a prediction for $\E_{c \sim \mathcal{D}_x}[c|x]$, which we denote by $\hat{c}$, and solving the deterministic version of the optimization problem based on $\hat{c}$, i.e., $ \min_{w \in S} \hat{c}^\top w $. Our primary interests in this paper concern defining suitable loss functions for the predict-then-optimize framework, examining their properties, and developing algorithms for training prediction models using these loss functions.



We now formally list the key ingredients of our framework:
\begin{enumerate}[leftmargin=*]
\item 
\emph{Nominal (downstream) optimization problem}, which is of the form
\begin{equation}\label{poi}
\begin{array}{rrccl}
P(c):~ & z^\ast(c) := & & \min\limits_{w} & c^Tw \\
& & & \text{s.t.} & w \in S \ ,
\end{array}
\end{equation}
where $w \in \bbR^d$ are the decision variables, $c \in \bbR^d$ is the problem data describing the linear objective function, and $S \subseteq \bbR^d$ is a nonempty, compact (i.e., closed and bounded), and convex set representing the feasible region. 
\edit{Since we are focusing on \emph{linear} optimization problems herein, the assumptions that $S$ is convex and closed are without loss of generality. Indeed, if $S$ in \eqref{poi} is instead possibly non-convex or non-closed, then replacing $S$ by its closed convex hull does not change the optimal value $z^\ast(c)$ (Lemma 8 in \cite{jaggi2011convex}). Thus, this basic equivalence for linear optimization problems implies that our methodology can be applied to combinatorial and mixed-integer optimization problems, which we elaborate further on in Section \ref{sec:surrogate}.
}
Since $S$ is assumed to be fixed and known with certainty, every problem instance can be described by the corresponding cost vector, hence the dependence on $c$ in \eqref{poi}. When solving a particular instance where $c$ is unknown, a prediction for $c$ is used instead. We assume access to a practically efficient optimization oracle, $w^*(c)$, that returns a solution of $P(c)$ for any input cost vector. For instance, if \eqref{poi} corresponds to a linear, conic, or a mixed-integer optimization problem, then a commercial optimization solver or a specialized algorithm suffices for $w^\ast(c)$.  
\item \emph{Training data} of the form $(x_1, c_1), (x_2, c_2), \ldots, (x_n, c_n)$, where $x_i \in \calX$ is a feature vector representing contextual information associated with $c_i$.  
\item A \emph{hypothesis class} $\calH$ of cost vector prediction models $f: \calX \to \bbR^d$, where $\hat c := f(x)$ is interpreted as the predicted cost vector associated with feature vector $x$. 
\item A \emph{loss function} $\ell(\cdot, \cdot) : \bbR^d \times \bbR^d \to \bbR_{+}$, whereby $\ell(\hat c, c)$ quantifies the error in making prediction $\hat c$ when the realized (true) cost vector is actually $c$.
\end{enumerate}

Given the loss function $\ell(\cdot, \cdot)$ and the training data $(x_1, c_1), \ldots, (x_n, c_n)$, the empirical risk minimization (ERM) principle states that we should determine a prediction model $f^\ast \in \calH$ by solving the optimization problem
\begin{equation}\label{erm1}
\min_{f \in \calH} \ \frac{1}{n}\sum_{i = 1}^n \ell(f(x_i), c_i) \ . 
\end{equation}
Provided with the prediction model $f^*$ and given a feature vector $x$, the predict-then-optimize decision rule is to choose the optimal solution with respect to the predicted cost vector, i.e., $w^\ast(f^*(x))$. Example \ref{nf} in Appendix \ref{sec:examples} contextualizes our framework in the context of a network optimization problem. 


In standard applications of the ``Predict, then Optimize'' framework, as in Example \ref{nf}, the loss function that is used is completely independent of the nominal optimization problem. In other words, the underlying structure of the optimization problem $P(\cdot)$ does not factor into the loss function and therefore the training of the prediction model. For example, when $\ell(\hat c, c) = \tfrac{1}{2}\|\hat c - c\|_2^2$, this corresponds to the least squares (LS) loss function. Moreover, if $\calH$ is a set of linear predictors, then \eqref{erm1} reduces to a standard least squares linear regression problem. In contrast, our focus in Section \ref{sec:spo} is on the construction of loss functions that measure decision errors in predicting cost vectors by leveraging problem structure.

\paragraph{Useful Notation.} Let $p$ be the dimension of a feature vector, $d$ be the dimension of a decision vector, and $n$ be the number of training samples. Let $W^\ast(c) := \arg\min_{w \in S}\left\{c^Tw\right\}$ denote the set of optimal solutions of $P(\cdot)$, and let $w^\ast(\cdot) : \bbR^d \to S$ denote a particular \emph{oracle} for solving $P(\cdot)$. That is, $w^\ast(\cdot)$ is a fixed deterministic mapping such that $w^\ast(c) \in W^\ast(c)$. Note that nothing special is assumed about the mapping $w^\ast(\cdot)$, hence $w^\ast(c)$ may be regarded as an arbitrary element of $W^\ast(c)$. 
Let $\supp(\cdot) : \bbR^d \to \bbR$ denote the support function of $S$, which is defined by $\supp(c) := \max_{w \in S}\{c^Tw\}$. Since $S$ is compact, $\supp(\cdot)$ is finite everywhere, the maximum in the definition is attained for every $c \in \bbR^d$, and note that $\supp(c) = -z^\ast(-c) = c^Tw^\ast(-c)$ for all $c \in \bbR^d$. Recall also that $\supp(\cdot)$ is a convex function. For a given convex function $h(\cdot) : \bbR^d \to \bbR$, recall that $g \in \bbR^d$ is a subgradient of $h(\cdot)$ at $c \in \bbR^d$ if $h(c^\prime) \geq h(c) + g^T(c^\prime - c)$ for all $c^\prime \in \bbR^d$, and the set of subgradients of $h(\cdot)$ at $c$ is denoted by $\partial h(c)$. For two matrices $B_1, B_2 \in \bbR^{d \times p}$, the trace inner product is denoted by $B_1 \bullet B_2 := \trace(B_1^TB_2)$. Finally, we note that the name of the framework is inspired by \cite{farias2007revenue}.


\section{SPO Loss Functions} \label{sec:spo}

\begin{figure}[h!]
\FIGURE{
    \centering
    \begin{subfigure}[t]{0.4\textwidth}
        \centering
        \includegraphics[scale=.5]{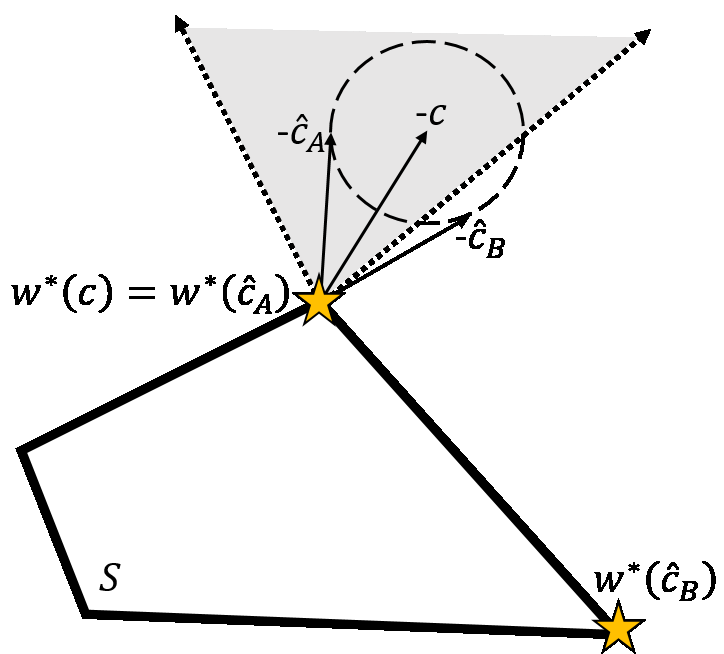}
        \caption{Polyhedral feasible region}
    \end{subfigure}
    \begin{subfigure}[t]{0.4\textwidth}
        \centering
        \includegraphics[scale=.5]{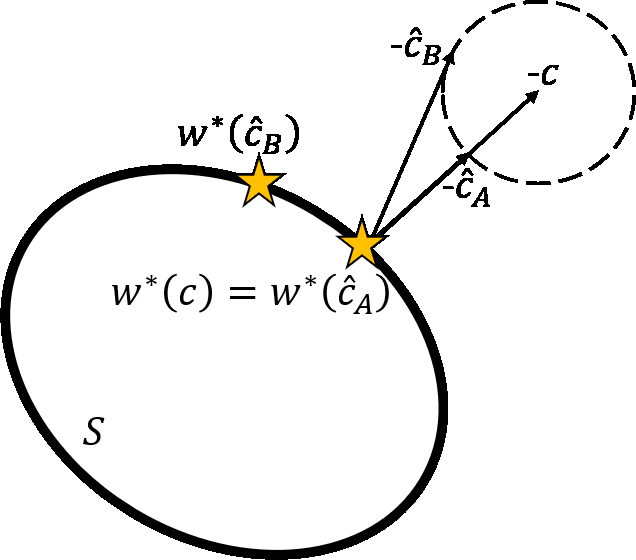}
        \caption{Elliptic feasible region}
    \end{subfigure}
}
{\edit{Geometric Illustration of SPO Loss\vspace{.3cm}}}
{\edit{In these two figures, we consider a two-dimensional polyhedron and ellipse for the feasible region $S$. We plot the (negative) of the true cost vector $c$, as well as two candidate predictions $\hat{c}_A$ and $\hat{c}_B$ that are equidistant from $c$ and thus have equivalent LS loss. One can see that the optimal decision for $\hat{c}_A$ coincides with that of $c$, since $w^*(\hat{c}_A)=w^*(c)$. In the polyhedron example, any predicted cost vector whose negative is not in the gray region will result in a wrong decision, where as in the ellipse example any predicted cost vector that is not exactly parallel with $c$ results in a wrong decision.}    \label{fig:geometry}}
\end{figure}

Herein, we introduce several loss functions that fall into the predict-then-optimize paradigm, but that are also \textit{smart} in that they take the nominal optimization problem $P(\cdot)$ into account when measuring errors in predictions. We refer to these loss functions as Smart ``Predict, then Optimize'' (SPO) loss functions. As a starting point, let us consider a true SPO loss function that exactly measures the excess cost incurred when making a suboptimal decision due to an imprecise cost vector prediction. Following the PO paradigm, given a cost vector prediction $\hat c$, a decision $w^\ast(\hat c)$ is implemented based on solving $P(\hat{c})$. After the decision $w^\ast(\hat c)$ is implemented, the cost incurred is with respect to the cost vector $c$ that is \emph{actually realized}. The excess cost due to the fact that $w^\ast(\hat c)$ may be suboptimal with respect to $c$ is then $c^Tw^\ast(\hat c) - z^\ast(c)$, which we call the SPO loss. \edit{In Figure \ref{fig:geometry}, we show how two predicted values of $c$ with the same prediction error can result in different decisions and different SPO losses. In fact, Figure \ref{fig:geometry} shows that the SPO loss can be 0 when $S$ is a polyhedron if $-\hat{c}$ lies in the cone corresponding to the extreme point $w^*(c)$, or when $S$ is an ellipse and $\hat{c}$ is in the same direction and parallel to $c$.} Definition \ref{true_def} formalizes this true SPO loss associated with making the prediction $\hat c$ when the actual cost vector is $c$, given a particular oracle $w^\ast(\cdot)$ for $P(\cdot)$.

\begin{definition}[SPO Loss] \label{true_def}
Given a cost vector prediction $\hat c$ and a realized cost vector $c$, the \emph{true SPO loss} $\ellTw(\hat c, c)$ w.r.t. optimization oracle $w^\ast(\cdot)$ is defined as $\ellTw(\hat c, c) :=~ c^Tw^\ast(\hat c) - z^\ast(c) $ .
\end{definition}

Note that there is an unfortunate deficiency in Definition \ref{true_def}, which is the dependence on the particular oracle $w^\ast(\cdot)$ used to solve \eqref{poi}. Practically speaking, this deficiency is not a major issue since we should usually expect $w^\ast(\hat c)$ to be a unique optimal solution, i.e., we should expect $W^*(\hat c)$ to be a singleton. Note that if any solution from $W^*(\hat{c})$ may be used by the loss function, then the loss function essentially becomes $\min_{w \in W^*(\hat{c})} c^T w -z^*(c)$. Thus, a prediction model would then be incentivized to always make the degenerate prediction $\hat{c}=0$ since $W^*(0)=S$. This would then imply that the SPO loss is 0.


In any case, if one wishes to address the dependence on the particular oracle $w^\ast(\cdot)$ in Definition \ref{true_def}, then it is most natural to ``break ties'' by presuming that the implemented decision has worst-case behavior with respect to $c$. Definition \ref{true_def2} is an alternative SPO loss function that does not depend on the particular choice of the optimization oracle $w^\ast(\cdot)$. 

\begin{definition}[Unambiguous SPO Loss]\label{true_def2}
Given a cost vector prediction $\hat c$ and a realized cost vector $c$, the (unambiguous) \emph{true SPO loss} $\ellT(\hat c, c)$ is defined as $
\ellT(\hat c, c) := \max_{w \in W^\ast(\hat c)}\left\{c^Tw\right\} - z^\ast(c) $.
\end{definition}

Note that Definition \ref{true_def2} presents a version of the true SPO loss that upper bounds the version from Definition \ref{true_def}, i.e., it holds that $\ellTw(\hat c, c) \leq \ellT(\hat c, c)$ for all $\hat c, c \in \bbR^d$. As mentioned previously, the distinction between Definitions \ref{true_def} and \ref{true_def2} is only relevant in degenerate cases. In the results and discussion herein, we work with the unambiguous true SPO loss given by Definition \ref{true_def2}. Related results may often be inferred for the version of the true SPO loss given by Definition \ref{true_def} by recalling that Definition \ref{true_def2} upper bounds Definition \ref{true_def} and that the two loss functions are almost always equal except for degenerate cases where $W^*(\hat{c})$ has multiple optimal solutions. 

Notice that $\ellT(\hat c, c)$ is impervious to the scaling of $\hat c$, in other words it holds that $\ellT(\alpha\hat c, c) = \ellT(\hat c, c)$ for all $\alpha > 0$. This property is intuitive since the true loss associated with prediction $\hat c$ should only depend on the optimal \emph{solution} of $P(\cdot)$, which does not depend on the scaling of $\hat c$. Moreover, this property is also shared by the 0-1 loss function in binary classification problems. Namely, labels can take values in the set $\{-1, +1\}$ and the prediction model predicts values in $\mathbb{R}$. If the predicted value has the same sign as the true value, the loss is 0, and otherwise the loss is 1. \edit{That is, given a predicted value $\hat c \in \bbR$ and a label $c \in \{-1, +1\}$, the 0-1 loss function is defined by $\ell_{0-1}(\hat c, c) := \mathbf{1}(\text{sgn}(\hat c) = c)$ where $\text{sgn}(\cdot)$ is the sign function and $\mathbf{1}(\cdot)$ is an indicator function equal to 1 if its input is true and 0 otherwise.}
Therefore, the 0-1 loss function is also independent of the scale on the predictions.  This similarity is not a coincidence; in fact, Proposition \ref{binary1} illustrates that binary classification is a special case of the SPO framework. All proofs can be found in Appendix \ref{sec:proofs}.

\begin{proposition}[SPO Loss Generalizes 0-1 loss]\label{binary1}
When $S = [-1/2, +1/2]$ and $c \in \{-1, +1\}$, then $\ellT(\hat{c}, c) = \mathbf{1}(\text{sgn}(\hat c) = c)$, i.e., the SPO loss function exactly matches  the 0-1 loss function associated with binary classification.
\end{proposition}

Now, given the training data, we are interested in determining a cost vector prediction model with minimal true SPO loss. Therefore, given the previous definition of the true SPO loss $\ellT(\cdot, \cdot)$, the prediction model would be determined by following the empirical risk minimization principle as in \eqref{erm1}, which leads to the following optimization problem:
\begin{equation}\label{erm_true}
\min_{f \in \calH} \ \frac{1}{n}\sum_{i = 1}^n \ellT(f(x_i), c_i) \ . 
\end{equation}
Unfortunately, the above optimization problem is difficult to solve, both in theory and in practice. Indeed, for a fixed $c$, $\ellT(\cdot, c)$ may not even be continuous in $\hat c$ since $w^\ast(\hat c)$ (and the entire set $W^\ast(\hat c)$) may not be continuous in $\hat c$. Moreover, since Proposition \ref{binary1} demonstrates that our framework captures binary classification, solving \eqref{erm_true} is at least as difficult as optimizing the 0-1 loss function, which may be NP-hard in many cases \citep{ben2003difficulty}. We are therefore motivated to develop approaches for producing ``reasonable'' approximate solutions to \eqref{erm_true} that {\em (i)} outperform standard PO approaches, and {\em (ii)} are applicable to large-scale problems where the number of training samples $n$ and/or the dimension of the hypothesis class $\cal H$ may be very large.


\subsection{An Illustrative Example}\label{sect:illustrative}

In order to build intuition, we now compare the SPO loss against the classical least squares (LS) loss function via an illustrative example. Consider a very simple shortest path problem with two nodes $s$ and $t$. There are two edge that go from $s$ to $t$, edge 1 and edge 2. Thus, a cost vector $c$ is 2-dimensional in this setting, and the goal is to simply choose the edge with the lower cost. We shall not observe $c$ directly at the decision-making time, but rather just a 1-dimensional feature $x$ associated with the vector $c$. Our data consists of $(x_i,c_i)$ pairs, and $c_i$ are generated nonlinearly as a function of $x_i$. 

\begin{figure}[h!]
\FIGURE{
    \centering
    \begin{subfigure}[t]{0.45\textwidth}
        \centering
        \includegraphics[width = \textwidth]{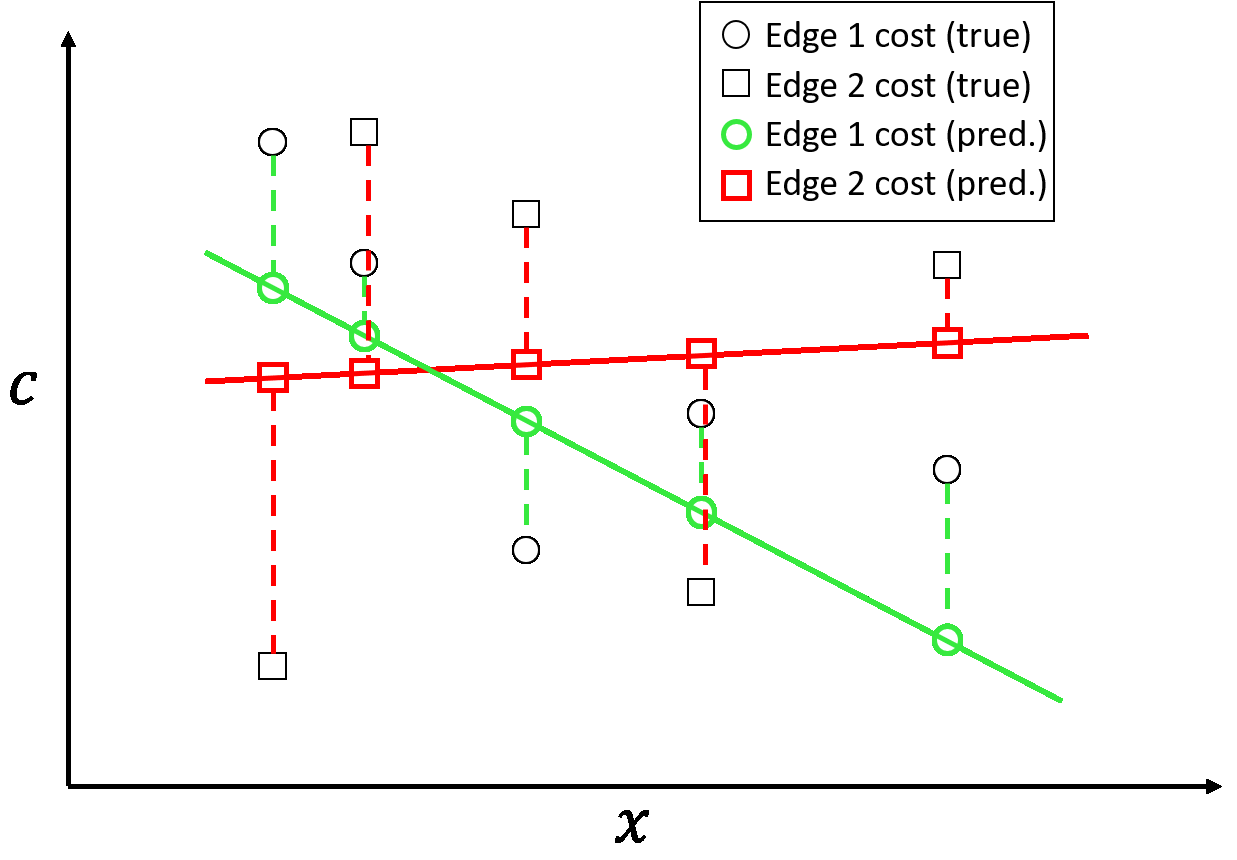}
        \caption{Prediction residuals}
    \end{subfigure}
    \begin{subfigure}[t]{0.45\textwidth}
        \centering
        \includegraphics[width = \textwidth]{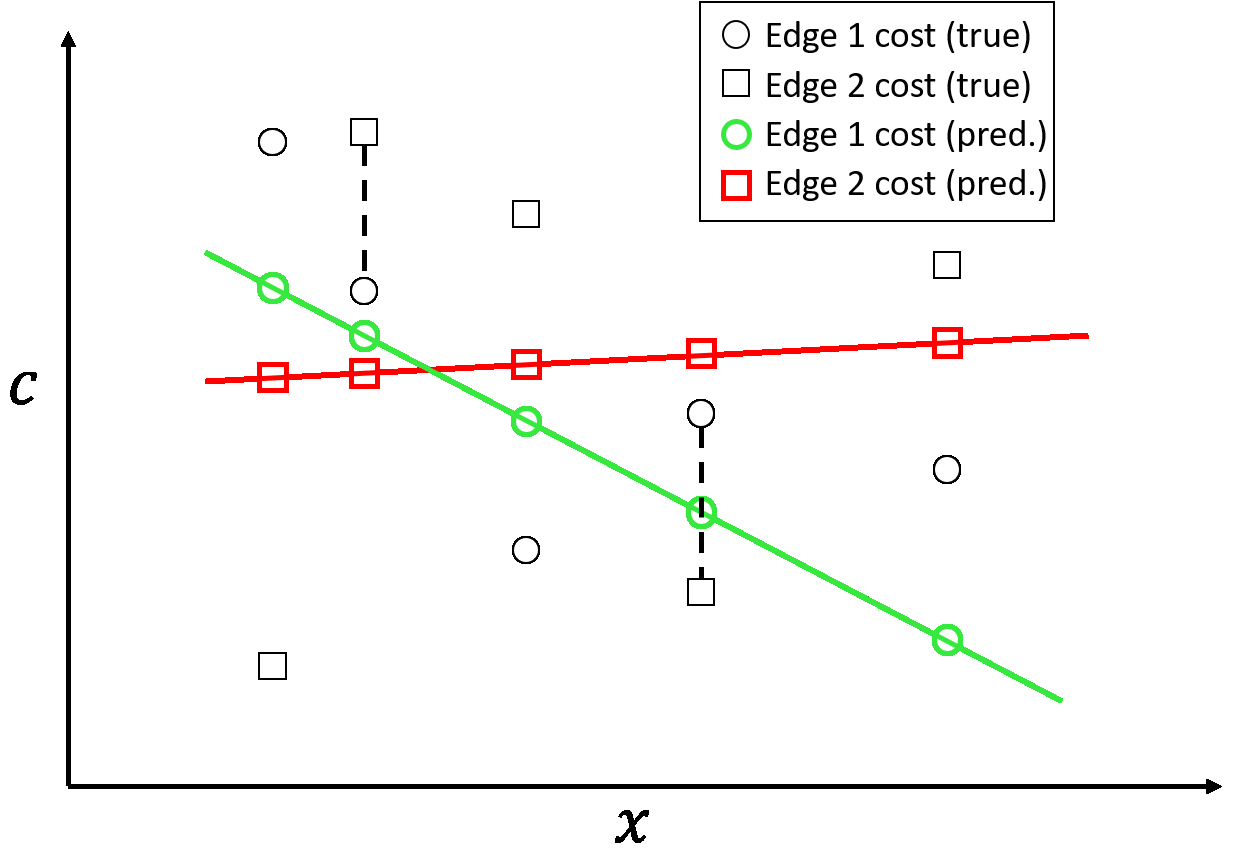}
        \caption{Decision residuals}
    \end{subfigure}
}
{\edit{Difference between prediction and decision residuals}}
{\edit{In (a), the residuals for the LS loss function are marked by the dashed lines. The residual is the distance between the prediction and the true value. In (b), the residuals for the SPO loss function are marked by the dashed black lines. The residual is 0 when the predicted values are in the right order. Otherwise, the residual is the distance between the true values.}
    \label{fig:residuals}}
\end{figure}

 The goal of the decision maker is to predict the cost of each edge from the feature using a simple linear regression model. The intersection of the two lines (corresponding to each edge) will signal the decision boundary in the predict-then-optimize framework. The decision maker shall try both the SPO and LS loss functions to do the linear regression. \edit{In Figure \ref{fig:residuals}, we illustrate the difference between LS and SPO by visualizing the residuals for one particular dataset and linear models for prediction the edge 1 and edge 2 costs. In LS regression, one minimizes the sum of the residuals squared, which is denoted by the dashed green and red lines in  Figure \ref{fig:residuals}(a). When using SPO loss, we consider ``decision residuals'' which only occur when the predictions result in choosing the wrong edge. In these cases, the SPO cost is the magnitude difference between the two true costs of edge 1 and edge 2, as depicted in Figure \ref{fig:residuals}(b).}
 
 \begin{figure}[h!]
\FIGURE
{
\centering
\includegraphics[scale=.5]{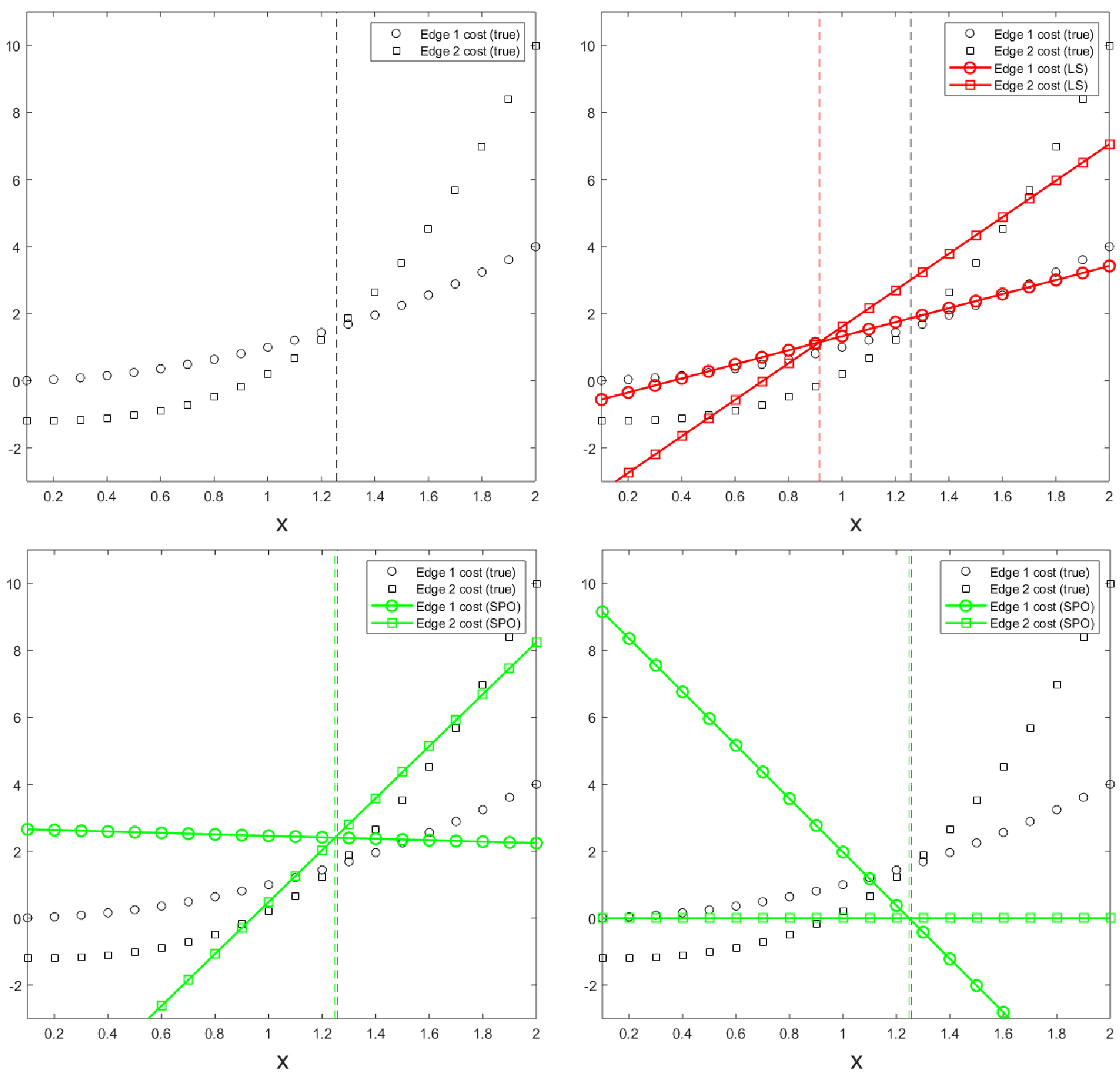}
}
{Illustrative Example. \vspace{.5cm} \label{fig:example}}
{The circles correspond to edge 1 costs and the squares correspond to edge 2 costs. Red lines and points correspond to the least squares fit and predictions, while green lines and points correspond to the SPO fit and predictions. The vertical dotted lines correspond to the decision boundaries under the true and prediction models. \edit{The SPO+ decision boundary in this stylized example coincides with the SPO decision boundary.}}
\end{figure}

 \edit{In Figure \ref{fig:example}, we consider another dataset, but this time plot the optimal LS and SPO linear regression models. In the first panel of Figure \ref{fig:example}, we plot the dataset and the optimal decision boundary. In the second panel, we plot the best LS fit to the data, and in the last two panels we plot two different optimal solutions to the SPO linear regression. (\edit{In fact, the SPO fitted models are also optimal for SPO+ loss which we derive in Section \ref{sec:surrogate}.}) Note that the the SPO loss in Figure \ref{fig:example} is 0, as there are no decision errors as described in Figure \ref{fig:residuals}.}
 
 One can see from Figure \ref{fig:example} that the LS lines very closely approximate the nonlinear data, although the decision boundary for LS is quite far from the optimal decision boundary. For any value of $x$ between the dotted black and red lines, the decision maker will choose the wrong edge. In contrast, the SPO lines need not approximate the data well at all, yet its decision boundary is nearly-optimal. In fact, the SPO lines have 0 training error, despite not fitting the data at all. \edit{The key intuition is that the SPO loss is incurred anytime the wrong edge is chosen, and in this example one can construct lines that cross at the right decision boundary so that the wrong edge is never chosen, resulting in zero SPO loss. Note that the only important consideration is where the lines intersect, and thus the SPO linear regression does not necessarily minimize prediction error.} Of course, a convex combination of SPO and LS loss may be used to overcome the unusual looking lines generated. In fact, there are infinitely optimal solutions to the ERM problem for the SPO loss, all of which just require that the intersection of the lines occurs between the $x$ values of 0.8 and 0.9.

\subsection{The SPO+ Loss Function} \label{sec:surrogate}

In this section, we focus on deriving a tractable surrogate loss function that reasonably approximates $\ellT(\cdot,\cdot)$. Our surrogate function $\ellS(\cdot,\cdot)$, which we call the SPO+ loss function, can be derived in a few steps that we shall carefully justify below. Ideally, when finding the prediction model that minimizes the empirical risk using the SPO+ loss, this prediction model will also approximately minimize \eqref{erm_true}, the empirical risk using the SPO loss.

To begin the derivation of the SPO+ loss, we first observe that for any $\alpha \in \mathbb{R}$, the SPO loss can be written as
\begin{align}
\ellT(\hat{c},c) = \max_{w \in W^\ast(\hat c)}\left\{c^Tw - \alpha \hat{c}^T w \right\} + \alpha z^*(\hat{c}) - z^\ast(c) \  \label{der1}
\end{align}
since $z^*(\hat{c})=\hat{c}^T w$ for all $w \in W^*(\hat{c})$. Clearly, replacing the constraint $w \in W^\ast(\hat c)$ with $w \in S$ in \eqref{der1} results in an upper bound. Since this is true for all values of $\alpha$, then 
\begin{align}
\ellT(\hat{c},c) ~\leq~ \inf_{\alpha} \left\{ \max_{w \in S}\left\{c^T w - \alpha \hat{c}^T w \right\} + \alpha z^*(\hat{c})  \right\} - z^\ast(c) \ . \label{der2}
\end{align}
In fact, one can show that inequality \eqref{der2} is actually an equality using duality theory, and moreover, the optimal value of $\alpha$ tends to $\infty$. Intuitively, one can see that as $\alpha$ gets large, then the term $c^T w$ in the inner maximization objective becomes negligible and the solution tends to $w^*(\alpha \hat{c})=w^*(\hat{c})$. Thus, as $\alpha$ tends to $\infty$, the inner maximization over $S$ can be replaced with maximization over $W^*(\hat{c})$, which recovers \eqref{der1}. We formalize this equivalence in Proposition \ref{alternative-rep} below. 

\begin{proposition}[Dual Representation of SPO Loss]\label{alternative-rep} 
For any cost vector prediction $\hat c \in \bbR^d$ and realized cost vector $c \in \bbR^d$, the function $\alpha \mapsto \max_{w \in S}\left\{c^T w - \alpha \hat{c}^T w \right\} + \alpha z^*(\hat{c})$ is monotone decreasing on $\bbR$, and the true SPO loss function  may be expressed as
\begin{equation}\label{represent_alpha_inf}
\ellT(\hat c, c) ~=~ \lim_{\alpha \rightarrow \infty} \left\{ \max_{w \in S}\left\{c^T w - \alpha \hat{c}^T w \right\} + \alpha z^*(\hat{c})  \right\} - z^\ast(c) \ .
\end{equation}
\end{proposition}
 
Using Proposition \ref{alternative-rep}, we shall now revist the SPO ERM problem \eqref{erm_true}  which can be written as 
\begin{align}
&~ \min_{f \in \calH} \ \frac{1}{n}\sum_{i = 1}^n \lim_{\alpha_i \rightarrow \infty} \left\{ \max_{w\in S} \left\{ c_i^T w - \alpha_i f(x_i)^T w \right\} + \alpha_i z^\ast(f(x_i)) \right\} - z^\ast(c_i)  \nonumber \\
=&~ \min_{f \in \calH} \ \frac{1}{n}\sum_{i = 1}^n \lim_{\alpha_i \rightarrow \infty} \left\{ \max_{w\in S} \left\{ c_i^T w - \alpha_i f(x_i)^T w \right\} + \alpha_i f(x_i)^T w^\ast(\alpha_i f(x_i)) \right\} - z^\ast(c_i) \nonumber \\
=&~ \min_{f \in \calH} \ \frac{1}{n} \lim_{\alpha \rightarrow \infty} \left\{ \sum_{i = 1}^n  \max_{w\in S} \left\{ c_i^T w - \alpha f(x_i)^T w \right\} + \alpha f(x_i)^T w^\ast(\alpha f(x_i))- z^\ast(c_i) \right\}  \nonumber \\ 
\leq&~ \min_{f \in \calH} \ \frac{1}{n}  \sum_{i = 1}^n  \max_{w\in S} \left\{ c_i^T w - 2 f(x_i)^T w \right\} + 2f(x_i)^T w^\ast(2 f(x_i))- z^\ast(c_i)  \label{der3}\\
\leq&~ \min_{f \in \calH} \ \frac{1}{n}  \sum_{i = 1}^n  \max_{w\in S} \left\{ c_i^T w - 2 f(x_i)^T w \right\} + 2f(x_i)^T w^\ast(c_i)- z^\ast(c_i) \ . \label{der4}
\end{align}
The first equality follows from the fact that $z^\ast(\alpha_if(x_i)) = \alpha_iz^\ast(f(x_i))$ for any $\alpha_i > 0$.
The second equality follows from the observation that all of the $\alpha_i$ variables are tending to the same value, so we can replace them with one variable which we call $\alpha$.  The first inequality follows from Proposition \ref{alternative-rep}, in particular that setting $\alpha=2$ in \eqref{der2} results in an upper bound on the SPO loss (we shall revisit this specific choice below). Finally, the second inequality follows from the fact that $w^*(c_i)$ is a feasible solution of $P(2 f(x_i))$.

The summand expression in \eqref{der4} is exactly what we refer to as the SPO+ loss function, which we formally state in Definition \ref{spo_def}.  
 
\begin{definition}[SPO+ Loss]\label{spo_def}
Given a cost vector prediction $\hat c$ and a realized cost vector $c$, the \emph{SPO+ loss} is defined as 
$\ellS(\hat c, c) := \max_{w \in S}\left\{c^T w - 2\hat{c}^T w \right\} + 2\hat{c}^Tw^\ast(c) - z^\ast(c)$.
\end{definition}

Recall that $\supp(\cdot)$ is the support function of $S$, i.e., $\supp(c) := \max_{w \in S}\{c^Tw\}$. Using this notation, the SPO+ loss may be equivalently expressed as $\ellS(\hat c, c) = \supp(c - 2\hat c) + 2\hat{c}^Tw^\ast(c) - z^\ast(c)$.

Before proceeding, we shall provide reasoning as to why inequalities \eqref{der3} and \eqref{der4}, which were used to derive SPO+, are indeed reasonable approximations. Although inequality \eqref{der3} could have been derived without the intermediary steps before it, we now claim that this inequality is actually an equality for many hypothesis classes. Namely, for any hypothesis class $\calH$ where $f\in \calH$ implies $\alpha f \in \calH$ for all $\alpha \geq 0$, then the inequality is tight since minimizing over $\alpha f$ is equivalent to minimizing over $2 f$. For example, the hypothesis class of linear models satisfies this property since all scalar multiples of linear models are also linear.  Note that $\alpha$ being absorbed into the hypothesis class was possible because the $\alpha_i$ terms in each summand can be replaced by a single $\alpha$ since they all tend to infinity. \textit{We specifically choose $\alpha=2$ (rather than any other positive scalar) because the Bayes risk minimizer of the SPO+ loss  (under some conditions) is exactly $\E [c|x]$ rather than a multiple of $\E[c|x]$.} This notion will be formalized in Section \ref{sec:consistency}.
 
The final step, \eqref{der4}, in the derivation of our convex surrogate SPO+ loss function involves approximating the concave (nonconvex) function $z^\ast(\cdot)$ with a first-order expansion. Namely, we apply the bound $z^*(2f(x_i)) = 2z^\ast(f(x_i)) \leq 2f(x_i)^Tw^\ast(c_i)$, which can be viewed as a first-order approximation of $z^\ast(f(x_i))$ based on a supergradient computed at $c_i$ (i.e., it holds that $w^\ast(c_i) \in \partial z^\ast(c_i)$). Note that if $f(x_i)=c_i$, then $\ellT(f(x_i),c_i)=\ellS(f(x_i),c_i)=0$ which implies that when minimizing SPO+, intuitively we are trying to get $f(x_i)$ to be close to $c_i$. Therefore, one might expect $w^\ast(c_i)$ to be a near-optimal solution to $P(2f(x_i))$ and thus inequality \eqref{der4} would be a reasonable approximation.  In fact, Section \ref{sec:consistency} provides a consistency property under some assumptions that would suggest the prediction $f(x_i)$ is indeed reasonably close to the expected value of $c_i$ if the prediction model is trained on a sufficiently large dataset.

Next, we state the following proposition which formally shows that  the SPO+ loss is an upper bound on the SPO loss and it is function is convex in $\hat{c}$.  Note that while the SPO+ loss is convex in $\hat c$, in general it is not differentiable since $\supp(\cdot)$ is not generally differentiable. {\blockedit However, Proposition \ref{prop:convex_ub} also shows that $2(w^\ast(c) - w^\ast(2\hat c - c))$ is a subgradient of the SPO+ loss, which is utilized in developing computational approaches in Section \ref{sec:algorithms}.  

\begin{proposition}[SPO+ Loss Properties]\label{prop:convex_ub}
Given a fixed realized cost vector $c$, it holds that:
\begin{enumerate}
\item $\ellT(\hat c, c) ~\leq~ \ellS(\hat c, c) \ \text{ for all } \hat c \in \bbR^d$,
\item $\ellS(\hat{c}, c)$ is a convex function of the cost vector prediction $\hat c$, and
\item For any given $\hat c$, $2(w^\ast(c) - w^\ast(2\hat c - c))$ is a subgradient of $\ellS(\cdot)$ at $\hat c$, i.e., $2(w^\ast(c) - w^\ast(2\hat c - c)) \in \partial\ellS(\hat c, c)$.
\end{enumerate}
\end{proposition}
}

 The convexity of the SPO+ loss function is also shared by the hinge loss function, which is a convex upper bound for the 0-1 loss function. Recall that the hinge loss given a prediction $\hat{c}$ is $\max \{0, 1-\hat{c} \}$ if the true label is $1$ and $\max \{0, 1+\hat{c} \}$ if the true label is $-1$. More concisely, the hinge loss can be written as $\max \{0, 1-c\hat{c}\}$ where $c \in \{-1, +1\}$ is the true label. The hinge loss is central to the support vector machine (SVM) method, where it is used as a convex surrogate to minimize 0-1 loss. Recall that, in this setting of binary classification, the SPO loss exactly captures the 0-1 loss as formalized in Proposition \ref{binary1}. In the same setting, it turns out that the SPO+ loss is equal to the hinge loss evaluated at $2\hat{c}$, i.e., twice the predicted value, which is formalized below in Proposition \ref{binary2}. This mild discrepancy is due to our choice of $\alpha = 2$ in the above derivation of the SPO+ loss; the alternative choice of $\alpha = 1$ would yield the hinge loss exactly.  

\begin{proposition}[SPO+ Loss Generalizes Hinge Loss]\label{binary2}
Under the same conditions as Proposition \ref{binary1}, namely when $S = [-1/2, +1/2]$ and $c \in \{-1, +1\}$, it holds that $\ellS(\hat{c}, c) = \max\{0, 1-2c\hat{c}\}$, i.e., the SPO+ loss function is equivalent to the hinge loss function associated with binary classification.
\end{proposition}

\begin{remark}[Connection to structured prediction] It is worth pointing out that the previously described construction of the SPO+ loss bears some resemblance to the construction of the structured hinge loss (\citet{taskar2004max, taskar2005learning, tsochantaridis2005large,nowozin2011structured}) in structured support vector machines (SSVMs). Moreover, our problem setting expands upon that of structured prediction by utilizing the objective cost of the nominal optimization problem to naturally define the SPO loss function. That is, if we define $w_i^\ast := w^\ast(c_i)$, then the modified dataset $(x_1, w_1^\ast), (x_2, w_2^\ast), \ldots, (x_n, w_n^\ast)$ may be regarded as the training data of a structured prediction problem. However, this reduction throws away valuable information about the cost vectors $c_i$, whereas the SPO+ loss function naturally exploits this information and upper bounds the SPO loss. Hence, our framework (and the surrogate SPO+ loss function) may be viewed as a type of refinement of the SSVM problem (and the structured hinge loss) to settings where there is a natural cost structure. Note that both the SPO+ loss and the structured hinge loss recover the regular hinge loss of binary classification as a special case. 
The hinge loss satisfies a key consistency property with respect to the 0-1 loss \citep{steinwart2002support}, which justifies its use in practice. In Section \ref{sec:consistency} we show a similar consistency result for the SPO+ loss with respect to the SPO loss under some mild conditions.
On the other hand, the structured hinge loss is often \emph{inconsistent} (see, e.g., the discussion around equation (11) in \citet{zhang2004statistical}), although there have been results on characterizing properties of consistent loss function in multiclass classification and structured prediction  \citep{zhang2004statistical, tewari2007consistency, osokin2017structured}.
\Halmos \end{remark}

{\blockedit
\begin{remark}[When $P(\cdot)$ is a combinatorial or mixed-integer problem] \label{remark:mixed} 
As mentioned previously, the assumptions that $S$ is convex and closed are without loss of generality since one can simply replace a possibly non-convex or non-closed set with its closed convex hull in \eqref{poi} without changing the optimal value $z^\ast(c)$. To be more concrete, suppose that $\tilde{S} \subseteq \bbR^d$ is a bounded but possibly non-convex or non-closed set and that $S$ is the closed convex hull of $\tilde{S}$. Suppose further that the the oracle $w^\ast(\cdot)$ returns an optimal solution in $\tilde{S}$, i.e., $w^\ast(c) \in \arg\min_{w \in \tilde{S}}c^Tw \subseteq \arg\min_{w \in S}c^Tw$ for all $c \in \bbR^d$. For example, if $\tilde{S}$ represents the feasible region of a combinatorial or mixed-integer optimization problem, then the oracle would correspond to a practically efficient algorithm for this problem. Then, using the fact that linear optimization on $\tilde{S}$ is equivalent to linear optimization on $S$, it is easy to see that the SPO and SPO+ loss functions  defined with respect to $\tilde{S}$ exactly equal the corresponding loss functions defined with respect to $S$. Finally, using Proposition \ref{prop:convex_ub}, one can use the oracle $w^\ast(c) \in \arg\min_{w \in \tilde{S}}c^Tw$ to compute subgradients of the SPO+ loss function, which can be utilized in computational approaches as described in Section \ref{sec:algorithms}. \Halmos
\end{remark}
}

Applying the ERM principle as in \eqref{erm_true} to the SPO+ loss yields the following optimization problem for selecting the prediction model:
\begin{equation}\label{erm_spo}
\min_{f \in \calH} \ \frac{1}{n}\sum_{i = 1}^n \ellS(f(x_i), c_i) \ . 
\end{equation}
Much of the remainder of the paper describes results concerning problem \eqref{erm_spo}. In Section \ref{sec:consistency} we demonstrate the aforementioned Fisher consistency result, in Section \ref{sec:algorithms} we describe several computational approaches for solving problem \eqref{erm_spo}, and in Section \ref{sec:computation} we demonstrate that \eqref{erm_spo} often offers superior practical performance over standard PO approaches. Next, we provide a theoretically motivated justification for using the SPO+ loss.

\section{Consistency of the SPO+ Loss Function} \label{sec:consistency}
{\blockedit In this section, we prove a fundamental consistency property, known as \textit{Fisher consistency}, to describe when minimizing the SPO+ loss is equivalent to minimizing the SPO loss. The Fisher consistency of a surrogate loss function means that under full knowledge of the data distribution and no restriction on the hypothesis class, the function that minimizes the surrogate loss also minimizes the true loss \citep{lin2004note,zou2008new}. One may also say that the surrogate loss is calibrated with the true loss \citep{bartlett2006convexity}. Our result is analagous to the well-known consistency results of the hinge loss and logistic loss functions with respect to the 0-1 loss -- minimizing hinge and logistic loss under full knowledge also minimizes the 0-1 loss -- and provides theoretical motivation for their success in practice. 

More formally, we let $\calD$ denote the distribution of $(x,c)$, i.e., $(x, c) \sim \calD$, and consider the population  version of the true SPO risk (Bayes risk) minimization problem:
\begin{equation}\label{truerisk_true}
\min_{f} \ \mathbb{E}_{(x,c) \sim \calD} [\ellT(f(x), c)].
\end{equation}
and the population  version of the SPO+ risk minimization problem:
\begin{equation}\label{truerisk_spo}
\min_{f} \ \mathbb{E}_{(x,c) \sim \calD} [\ellS(f(x), c)] \ .
\end{equation}
Note here that we place no restrictions on $f(\cdot)$, meaning $\cal H$ consists of any measurable function mapping features to cost vectors. 
\begin{definition}[Fisher Consistency] A loss function $\ell(\cdot,\cdot)$ is said to be \textit{Fisher consistent} with respect to the SPO loss if $\argmin_{f}  \mathbb{E}_{(x,c) \sim \calD} [\ell(f(x), c)]$ (the set of minimizers of the Bayes risk of $\ell$) also minimizes \eqref{truerisk_true}.
\end{definition}

To gain some intuition, let $f^{\ast}_{\mathrm{SPO}}$ and $f^{\ast}_{\mathrm{SPO+}}$ denote any optimal solution of \eqref{truerisk_true} and \eqref{truerisk_spo}, respectively. From \eqref{core}, one can see that an ideal value for $f^*_{\mathrm{SPO}}(x)$ is simply $\E[c|x]$. In fact, as long as the optimal solution of $P(\E[c|x])$ is unique with probability 1 (over the distribution of $x \in \calX$), i.e., almost surely, then it is indeed the case $\bbE[c | x]$ is a minimizer of \eqref{truerisk_true} (see Proposition \ref{true_consistency} below). Moreover, any function that is almost surely equal to $\bbE[c | x]$ is also a minimizer of \eqref{truerisk_true}.
In Theorem \ref{thm:main}, we show that under Assumption \ref{assumption-fun}, any minimizer of the SPO+ population risk \eqref{truerisk_spo} must satisfy $f^{\ast}_{\mathrm{SPO+}}(x) = \E[c|x]$ almost surely and therefore also minimizes the SPO risk \eqref{truerisk_true}. In summary, the SPO+ loss is Fisher consistent with the SPO loss, under Assumption \ref{assumption-fun}.

\begin{assumption}\label{assumption-fun} These assumptions imply Fisher consistency of the SPO+ loss function.
\begin{enumerate}
\item Almost surely, $W^*(\E[c|x])$ is a singleton, i.e., $\mathbb{P}_x(|W^\ast(\E[c|x])|=1)=1$. 
\item For all $x\in \calX$, the distribution of $c|x$ is centrally symmetric about its mean $\E[c|x]$.
\item For all $x\in \calX$, the distribution of $c|x$ is continuous on all of $\mathbb{R}^d$.
\item The interior of the feasible region $S$ is nonempty.
\end{enumerate} 
\end{assumption}

\begin{theorem}[Fisher Consistency of SPO+]\label{thm:main}
Suppose Assumption \ref{assumption-fun} holds. Then, any minimizer of the SPO+ risk \eqref{truerisk_spo} is almost surely (over the distribution of $x \in \calX$) equal to $\E[c|x]$ and is also a minimizer of the SPO risk \eqref{truerisk_true}. Thus, the SPO+ loss function is Fisher consistent with respect to the SPO loss. 
\end{theorem}

The key results to prove Theorem \ref{thm:main} are provided in Section \ref{sec:proof}, and the final proof is given in the Appendix. We remark that Assumption \ref{assumption-fun}.1  is only needed to show that $\E[c|x]$ is a minimizer of the SPO risk. This assumption is rather mild as the set of points with multiple optimal solutions typically has measure 0. In fact, Assumption \ref{assumption-fun}.1  can be removed if one uses Definition \ref{true_def} of the SPO loss which uses a given optimization oracle. Assumption \ref{assumption-fun}.2 ensures that $\E[c|x]$ is a minimizer of the SPO+ risk. Note that a random vector $d$ is centrally symmetric about its mean if $d-\E[d]$ is equal in distribution to $\E[d]-d$, or equivalently $d$ is equal in distribution to $2\E[d] -d$. This symmetry condition is satisfied, for instance, when the data is assumed to be of the form $f(x) + \epsilon$ where $\epsilon$ is a zero-mean Gaussian distribution with a positive semi-definite covariance matrix. Finally, Assumptions \ref{assumption-fun}.3 and \ref{assumption-fun}.4, both of which are standard, are used to show that $\E[c|x]$ uniquely minimizes the SPO+ risk except possibly on a set of probability measure zero. Note that Assumptions \ref{assumption-fun}.2 and \ref{assumption-fun}.3 may be relaxed to hold almost surely with respect to the probability measure of $x \in \calX$; but for ease of presentation we state them for all $x \in \calX$. 
In Section \ref{sec:proof}, we discuss examples (provided in the Appendix) that show how our result may not hold if one of the assumptions are violated.

As mentioned previously, any minimizer for the least squares (LS) risk is also almost surely equal to $\E[c|x]$, and thus the least squares loss is also Fisher consistent with respect to the SPO loss. \textit{Thus, a priori, one cannot claim LS or SPO+ to be better than the other.} Indeed, we have derived a natural surrogate loss function, SPO+, directly from the SPO loss that maintains a fundamental consistency property of the de facto standard LS loss function. In fact, it is easy to see that under Assumption 1, any convex combination of the LS and SPO+ loss functions is Fisher consistent. Since this consistency property applies under full distributional information and no model misspecification (no restriction on hypothesis class), we show in Section \ref{sec:computation} that SPO+ indeed outperforms LS in several experimental settings, due to its ability to tailor the prediction to the optimization task.

}

\subsection{Key Results to Prove Fisher Consistency} \label{sec:proof}



{\blockedit
Throughout this section, we consider a non-parametric setup where the dependence on the features $x$ is dropped without loss of generality. To see this, first observe that the SPO risk satisfies $\mathbb{E}_{(x,c) \sim \calD} [\ellT(f(x), c)] = \bbE_x\left[\bbE_c\left[\ellT(f(x), c) ~|~ x\right]\right]$ and likewise for the SPO+ risk. Since there is no constraint on $f(\cdot)$ (the hypothesis class consists of all prediction models), solving problems \eqref{truerisk_true} and \eqref{truerisk_spo} is equivalent to optimizing each function value $f(x)$ individually for all $x \in \calX$. Therefore, for the remainder of the section unless otherwise noted, we drop the dependence on $x$. Thus, we now assume that the distribution $\calD$ is  only over $c$, and the SPO and SPO+ risk is defined as $\riskT(\hat c) := \bbE_{c} [\ellT(\hat c, c)]$ and $\riskS(\hat c) := \bbE_{c} [\ellS(\hat c, c)]$, respectively.  
For convenience, let us define $\bar{c} := \bbE_c[c]$ (note that we are implicitly assuming that $\bar c$ is finite).

Next, we fully characterize the minimizers of the true SPO risk problem \eqref{truerisk_true} in this setting.  Proposition \ref{true_consistency} demonstrates that for any minimizer $c^*$ of $\riskT(\cdot)$, all of its corresponding solutions with respect to the nominal problem, $W^*(c^*)$,  are also optimal solutions for $P(\bar c)$. In other words, minimizing the true SPO risk also optimizes for the expected cost in the nominal problem (since the objective function is linear). Proposition \ref{true_consistency} also demonstrates that the converse is true -- namely any cost vector prediction with a unique optimal solution that also optimizes for the expected cost is also a minimizer of the true SPO risk.

\begin{proposition}[SPO Minimizer] \label{true_consistency}
If a cost vector $c^\ast$ is a minimizer of $\riskT(\cdot)$, then $W^\ast(c^\ast) \subseteq W^\ast(\bar c)$. Conversely, if $c^\ast$ is a cost vector such that $W^\ast(c^\ast)$ is a singleton and $W^\ast(c^\ast) \subseteq W^\ast(\bar c)$, then $c^\ast$ is a minimizer of $\riskT(\cdot)$. 
\end{proposition}

 Example \ref{normal-example} in Appendix \ref{sec:examples} demonstrates that, in order to ensure that $c^\ast$ is a minimizer of $\riskT(\cdot)$, it is not sufficient to allow $c^\ast$ to be any cost vector such that $W^\ast(c^\ast) \subseteq W^\ast(\bar c)$. In fact, it may not be sufficient for $c^*$ to be $\bar{c}$. This follows from the unambiguity of the SPO loss function, which chooses a worst-case optimal solution in the event that the prediction allows for more than one optimal solution.

Next, we provide Proposition \ref{spo_consistent} which shows sufficient conditions for   $\bar{c}$ to be the minimizer of the SPO+ risk and therefore the minimizer of the SPO risk, implying Fisher consistency. We also provide conditions for when $\bar{c}$ is the unique minimizer of the SPO+ risk, which alleviates any concern that there may be alternate minimizers of the SPO+ risk which are not Fisher consistent.

\begin{proposition}[SPO+ Minimizer] \label{spo_consistent} Suppose that the distribution $\mathcal{D}$ of $c$ is continuous and centrally symmetric about its mean $\bar{c}$ (i.e., $c$ is equal in distribution to $2\bar{c} - c$). 
\begin{enumerate} 
\item[a)] Then $\bar{c}$ minimizes  $\riskS(\cdot)$.
\item[b)] In addition, suppose the interior of  $S$ is nonempty. Then  $\bar{c}$ is the unique minimizer of  $\riskS(\cdot)$.
\end{enumerate}
\end{proposition}

The two important assumptions in Proposition \ref{spo_consistent} are that $\mathcal{D}$ is centrally symmetric about its mean and continuous, both of which are not individually sufficient to ensure consistency on their own. Example \ref{bad_example2} in Appendix \ref{sec:examples}  demonstrates a situation where $c$ is continuous on $\mathbb{R}^d$ and the minimizer of SPO+ is unique, but it does not minimize the SPO risk. Example \ref{bad_example1} in Appendix \ref{sec:examples} demonstrates a situation where the distribution of $c$ is symmetric about its mean but there exists a minimizer of the SPO+ risk that does not minimize the SPO risk. Example \ref{interior} in Appendix \ref{sec:examples} demonstrates a case where the minimizer of SPO+ is not unique if $S$ is empty while $c$ is continuous and centrally symmetric about its mean. 

}
\section{Computational Approaches}\label{sec:algorithms}
In this section, we consider computational approaches for solving the SPO+ ERM problem \eqref{erm_spo}. Herein, we focus on the case of linear predictors, $\calH = \{f : f(x) = Bx \text{ for some } B \in \bbR^{d \times p}\}$, with regularization possibly incorporated into the objective function, using the regularizer $\Omega(\cdot) : \bbR^{d \times p} \to \bbR$. (This is equivalent to working with the hypothesis class $\calH = \{f : f(x) = Bx \text{ for some } B \in \bbR^{d \times p}, \Omega(B) \leq \rho \}$ for some $\rho > 0$.)
For example, we may use the ridge penalty $\Omega(B) = \tfrac{1}{2}\|B\|_F^2$, where $\|B\|_F$ denotes the Frobenius norm of $B$, i.e., the entry-wise $\ell_2$ norm. Other possibilities include an entry-wise $\ell_1$ penalty or the nuclear norm penalty, i.e., an $\ell_1$ penalty on the singular values of $B$. In any case, these presumptions lead to the following version of \eqref{erm_spo}:
\begin{equation}\label{erm_true3}
\min_{B \in \bbR^{d \times p}} \ \frac{1}{n}\sum_{i = 1}^n \ellS(Bx_i, c_i) ~+ \lambda\Omega(B) \ ,
\end{equation}
where $\lambda \geq 0$ is a regularization parameter. Since the SPO loss is convex as stated in Proposition \ref{prop:convex_ub}, then the above problem is a convex optimization problem as long as $\Omega(\cdot)$ is a convex function.

We mainly consider two approaches for solving problem \eqref{erm_true3}:  {\em (i)} reformulations based on modeling $\ellS(\cdot, c)$ using duality, and {\em (ii)} stochastic gradient based methods that instead rely only on an optimization oracle for problem \eqref{poi}. 
The reformulation based approach {\em (i)} requires an explicit description of the feasible region $S$, for example if $S$ is a polytope then this approach necessitates working with an explicit list of inequality constraints describing $S$. 
{\blockedit
On the other hand, the stochastic gradient based approach {\em (ii)} does not require an explicit description of $S$ and instead \emph{only} relies on iteratively calling the optimization oracle $w^\ast(\cdot)$ in order to compute stochastic subgradients of the SPO+ loss (see Proposition \ref{prop:convex_ub}). Therefore it is much more straightforward to apply the stochastic gradient descent approach to problems with complicated constraints, such as nonlinear problems as well as combinatorial and mixed-integer problems as mentioned in Remark \ref{remark:mixed}.}
While approach {\em (i)} is more restrictive in its requirements, it does offer 
a few advantages. Depending on the structure of $S$, for example if $S$ is a polytope with known linear inequality constraints, then approach {\em (i)} may able to utilize off-the-shelf conic optimization solvers such as CPLEX and Gurobi that are capable of producing high accuracy solutions for small to medium sized problem instances (see Section~\ref{sec:reform}).
However, for large scale instances where $d$, $p$, and $n$ might be very large, conic solvers based on interior point methods do not scale as well. Stochastic gradient methods, on the other hand, scale much better to instances where $n$ may be extremely large, and possibly also to instances where $d$ and $p$ are large but the optimization oracle $w^\ast(\cdot)$ is efficiently computable due to the special structure of $S$. The details of the approach \textit{(ii)} can be found in Appendix \ref{sec:sgd}.

\subsection{Reformulation Approach} \label{sec:reform}
We now discuss the reformulation approach {\em (i)}, which aims to recast problem \eqref{erm_true3} in a form that is amenable to popular optimization solvers. To describe this approach, we presume that $S$ is a polytope described by known linear inequalities, i.e., $S = \{w : Aw \geq b\}$ for some given problem data $A \in \bbR^{m \times d}$ and $b \in \bbR^m$. The same approach may also be applied to particular classes of nonlinear feasible regions, although the complexity of the resulting reformulated problem will be different.
The key idea is that when $S$ is a polytope, then $\ellS(\cdot, c)$ is a (piecewise linear) convex function of the prediction $\hat c$ and therefore the epigraph of $\ellS(\cdot, c)$ can be tractably modeled with linear constraints by employing linear programming duality.
Proposition \ref{reformulation} formalizes this approach. (Recall that, for $w \in \bbR^d$ and $x \in \bbR^p$, $wx^T$ denotes $d \times p$ outer product matrix where $(wx^T)_{ij} = w_ix_j$.)

\begin{proposition}[Reformulation of ERM for SPO+] \label{reformulation}
\edit{Suppose $S = \{w : Aw \geq b\}$ is a polytope.} Then the regularized SPO+ ERM problem \eqref{erm_true3} is equivalent to the following optimization problem:
\begin{equation}\label{reform}
\begin{array}{clr}
\min\limits_{B, p} & \displaystyle\frac{1}{n}\sum_{i = 1}^n \left[-b^Tp_i + 2(w^\ast(c_i)x_i^T) \bullet B - z^\ast(c_i)\right] ~+~ \lambda \Omega(B) \\ 
\mathrm{s.t.} & A^Tp_i = 2Bx_i - c_i & \mathrm{for \ all} \ i \in \{1, \ldots, n\} \\
& p_i \in \bbR^{m}, p_i \geq 0 & \mathrm{for \ all} \ i \in \{1, \ldots, n\} \\
& B \in \bbR^{d \times p} \ .
\end{array}
\end{equation}
\end{proposition}

\edit{Thus, as we can see, problem \eqref{reform} is almost a linear optimization problem -- the only part that may be nonlinear is the regularizer $\Omega(\cdot)$. For several natural choices of $\Omega(\cdot)$, problem \eqref{reformulation} may be cast as a conic optimization problem that can be solved efficiently with interior point methods. For instance, for the LASSO penalty where $\Omega(B) = \|B\|_1$, then \eqref{reform} is equivalent to a linear program.  If $\Omega(\cdot)$ is the ridge penalty, $\Omega(B) = \tfrac{1}{2}\|B\|_F^2$, then \eqref{reform} is equivalent to a quadratic program. If $\Omega(\cdot)$ is the nuclear norm penalty, $\Omega(B) = \|B\|_\ast$, then \eqref{reform} is equivalent to a semidefinite program.}



\section{Computational Experiments}\label{sec:computation}
In this section, we present computational results of synthetic data experiments wherein we empirically examine the quality of the SPO+ loss function for training prediction models, using the shortest path problem and portfolio optimization as our exemplary problem classes.
Following Section \ref{sec:algorithms}, we focus on linear prediction models, possibly with either ridge or entrywise $\ell_1$ regularization.
We compare the performance of four different methods:  
\begin{enumerate}
\item the previously described SPO+ method, \eqref{erm_true3}.
\item the least squares method that replaces the SPO+ loss function in \eqref{erm_true3} with $\ell(\hat c, c) = \tfrac{1}{2}\|\hat c - c\|_2^2$ and also uses regularization whenever SPO+ does.
\item an absolute loss function (i.e., $\ell_1$) approach that replaces the SPO+ loss function in \eqref{erm_true3} with $\ell(\hat c, c) = \|\hat c - c\|_1$ and also uses regularization whenever SPO+ does.
\item a random forests approach that independently trains $d$ different random forest models for each component of the cost vector, using standard parameter settings of $\lceil p/3 \rceil$ random features at each split and $100$ trees.
\end{enumerate}

Note that methods (2.), (3.) and (4.) above do not utilize the structure of $S$ in any way and hence may be viewed as independent learning algorithms with respect to each of the components of the cost vector. For methods (1.), (2.), and (3.) above, we include an intercept column in $B$ that is not regularized.
In order to ultimately measure and compare the performance of the four different methods, we compute a ``normalized'' version of the SPO loss of each of the four previously trained models on an independent test set of size $10,000$. Specifically, if $(\tilde x_1, \tilde c_1), (\tilde x_2, \tilde c_2), \ldots, (\tilde x_{n_{\text{test}}}, \tilde c_{n_{\text{test}}})$ denotes the test set, then we define the normalized test SPO loss of a previously trained model $\hat f$ by $\mathrm{NormSPOTest}(\hat f) := \frac{\sum_{i = 1}^{n_{\text{test}}} \ellT(\hat f(\tilde x_i), \tilde c_i)}{\sum_{i = 1}^{n_{\text{test}}} z^\ast(\tilde c_i)}$. Note that we naturally normalize by the total optimal cost of the test set given full information, which with high probability will be a positive number for the examples studied herein.

\subsection{Shortest Path Problem}
We consider a shortest path problem on a $5 \times 5$ grid network, where the goal is to go from the northwest corner to the southeast corner and the edges only go south or east. In this case, the feasible region $S$ can be modeled using network flow constraints as in Example \ref{nf}.
We utilize the reformulation approach given by Proposition \ref{reformulation} to solve the SPO+ training problem \eqref{erm_true3}. Specifically, we use the JuMP package in Julia (\cite{dunning2017jump}) with the Gurobi solver to implement problem \eqref{reform}. The optimization problems required in methods (2.) and (3.) are also solved directly using Gurobi. In some cases we use $\ell_1$ regularization for methods (1.), (2.), and (3.), in which case, in order to tune the regularization parameter $\lambda$, we try 10 different values of $\lambda$ evenly spaced on the logarithmic scale between $10^{-6}$ and $100$. Furthermore, we use a validation set approach where we train the $10$ different models on a training set of size $n$ and then use an independent validation set of size $n/4$ to pick the model that performs best with respect to the SPO loss.

\paragraph{Synthetic Data Generation Process.}
Let us now describe the process used for generating the synthetic experimental data instances for both problem classes.
Note that the dimension of the cost vector $d = 40$ corresponds to the total number of edges in the $5 \times 5$ grid network and that $p$ is a given number of features.
First, we generate a random matrix $B^\ast \in \bbR^{d \times p}$ that encodes the parameters of the true model, whereby each entry of $B^\ast$ is Bernoulli random variable that is equal to 1 with probability $0.5$. 
We generate the training data $(x_1, c_1), (x_2, c_2), \ldots, (x_n, c_n)$ and the testing data $(\tilde x_1, \tilde c_1), (\tilde x_2, \tilde c_2), \ldots, (\tilde x_n, \tilde c_n)$ according to the following generative model:
\begin{enumerate}
\item First, the feature vector $x_i \in \bbR^p$ is generated from a multivariate Gaussian distribution with i.i.d. standard normal entries, i.e., $x_i \sim N(0, I_p)$.
\item Then, the cost vector $c_i$ is generated according to $c_{ij} = \left[\left(\tfrac{1}{\sqrt{p}}(B^\ast x_i)_j + 3\right)^{\degg} + 1\right]\cdot\varepsilon_i^j$ for $j = 1, \ldots, d$, and where $c_{ij}$ denotes the $j^{\text{th}}$ component of $c_i$ and $(B^\ast x_i)_j$ denotes the $j^{\text{th}}$ component of $B^\ast x_i$. Here, $\degg$ is a fixed positive integer parameter and $\varepsilon_i^j$ is a multiplicative noise term that is generated independently at random from the uniform distribution on $[1 - \bar\varepsilon, 1 + \bar\varepsilon]$ for some parameter $\bar\varepsilon \geq 0$. 
\end{enumerate}

Note that the model for generating the cost vectors employs a polynomial kernel function (see, e.g., \cite{hofmann2008kernel}), whereby the regression function for the cost vector given the features, i.e., $\bbE[c|x]$, is a polynomial function of $x$ and the parameter $\degg$ dictates the degree of the polynomial. Importantly, we still employ a linear hypothesis class for methods (1.)-(3.) above, hence the parameter $\degg$ controls the amount of \emph{model misspecification} and as $\degg$ increases we expect the performance of the SPO+ approach to improve relative to methods (2.) and (3.). When $\degg = 1$, the expected value of $c$ is indeed linear in $x$. Furthermore, for large values of $\degg$, the least squares method will be sensitive to outliers in the cost vector generation process, which is our main motivation for also comparing against the absolute loss approach that is less sensitive to outliers. On the other hand, the random forests method is a non-parametric learning algorithm and will accurately learn the regression function for any value of $\degg$. However, the practical performance of random forests depends heavily on the sample size $n$ and, for relatively small values of $n$, random forests may perform poorly.

\begin{figure}[h!]
\centering
\includegraphics[width=\linewidth]{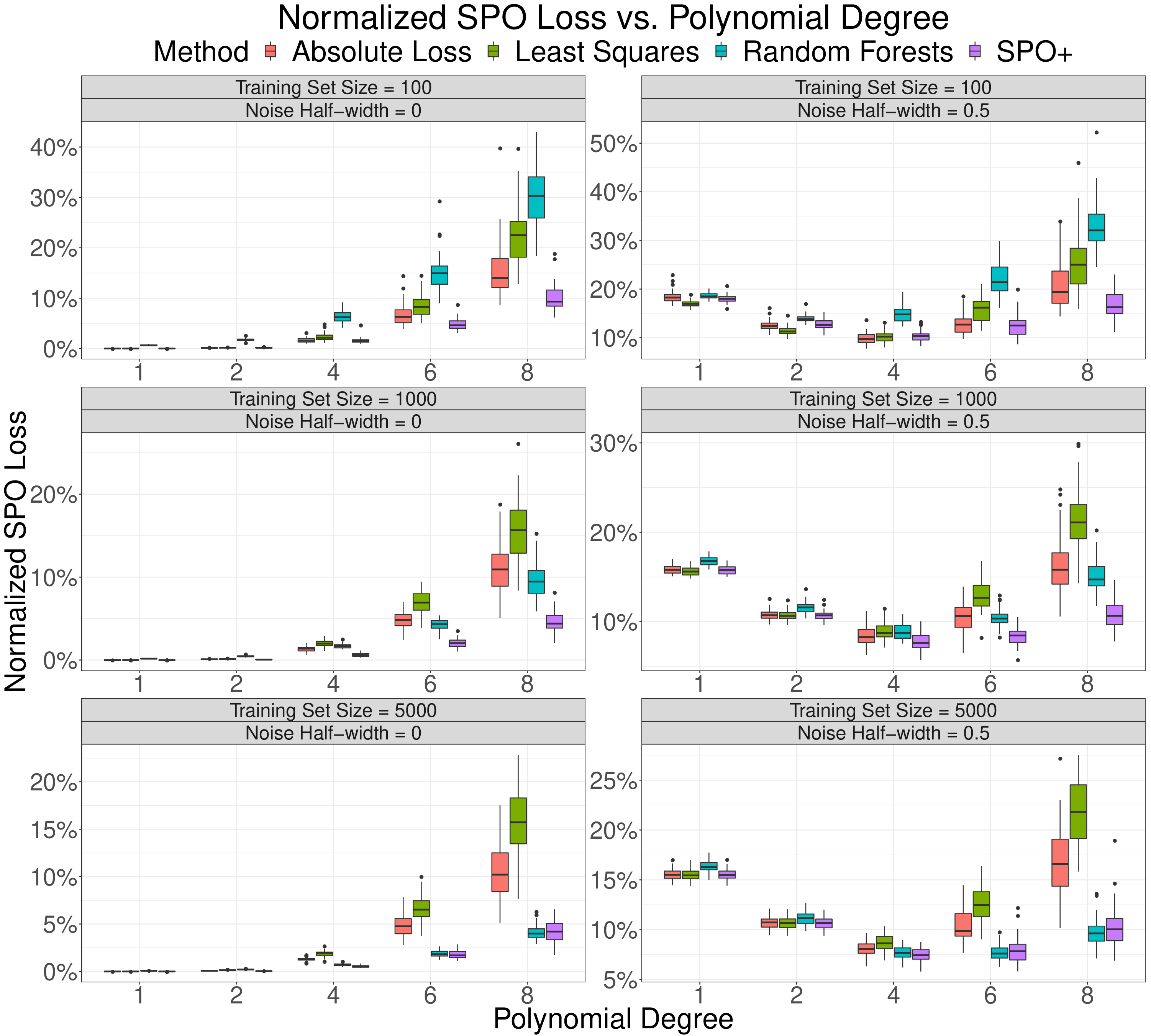}
\caption{Normalized test set SPO loss for the SPO+, least squares, absolute loss, and random forests methods on shortest path problem instances.}
\label{fig:adamexpts}
\end{figure}

\paragraph{Results.} In the following set of experiments on the shortest path problem we described, we fix the number of features at $p=5$ throughout and, as previously mentioned, use a $5 \times 5$ grid network, which implies that $d = 40$. Hence, in total there are $pd=200$ parameters to estimate. We vary the training set size $n \in \{100, 1000, 5000\}$, we vary the parameter $\degg \in \{1, 2, 4, 6, 8\}$, and we vary the noise half-width parameter $\bar{\varepsilon} \in \{0,  0.5\}$. For every value of $n$, $\degg$, and $\bar{\varepsilon}$, we run 50 simulations, each of which has a different $B^\ast$ and therefore different ground truth model. For the cases where $n \in \{100, 1000\}$, we employ $\ell_1$ regularization for methods (1.)-(3.), as previously described. When $n = 5000$ we do not use any regularization (since it did not appear to provide any value). As mentioned previously, for each simulation, we evaluate the performance of the trained models by computing the normalized SPO loss on a test set of 10,000 samples. \edit{The computation time for solving one ERM problem using the SPO+ loss is approximately 0.5-1.0 seconds, 5-30 seconds, and 1-15 minutes for $n \in \{100, 1000, 5000\}$, respectively. The other methods can be solved in a few seconds using well-developed packages.} Figure \ref{fig:adamexpts} summarizes our findings, and note that the box plot for each configuration of the parameters is across the 50 independent trials.


From Figure \ref{fig:adamexpts}, we can see that for small values of the $\degg$ parameter, i.e., $\degg \in \{1, 2\}$, the absolute loss, least squares, and SPO+ methods perform comparably, with the least squares method slightly dominating in the case of noise with $\bar{\varepsilon} = 0.5$. The slight dominance of least squares (and sometimes the absolute loss as well) in these cases might be explained by some inherent robustness properties of the least squares loss.
It is also plausible that, since the SPO+ loss function is more intricate than the ``simple'' least squares loss function, it may overfit in situations with noise and a small training set size.
On the other hand, as the parameter $\degg$ grows and the degree of model misspecification increases, then the SPO+ approach generally begins to perform best across all instances except when $n = 5000$, in which case random forests performs comparably to SPO+. This behavior suggests that the SPO+ loss is better than the competitors at leveraging additional data and stronger nonlinear signals.

It is interesting to point out that random forests generally does not perform well except when $n = 5000$, in which case it performs comparably to SPO+, which uses a much simpler linear hypothesis class. Indeed, when $n \in \{100, 1000\}$, random forests almost always performs worst, except for when $n = 1000$ and $\degg \in \{6, 8\}$, in which case random forests outperforms least squares, performs comparably to the absolute loss method, and is strongly dominated by SPO+. Indeed, the cases where $n \in \{1000, 5000\}$ and $\degg \in \{6, 8\}$ suggest that least squares is prone to outliers whereas the absolute loss is not, random forests is slow to converge due to its non-parametric nature, and SPO+ is best able to adapt to the large degree of model misspecification even with a modest amount of data (i.e., $n = 1000$).

\subsection{Portfolio Optimization}
Here we consider a simple portfolio selection problem based on the classical Markowitz model \citep{markowitz1952portfolio}.
As discussed in Section \ref{sec:intro}, we presume that there are auxiliary features that may be used to predict the returns of $d$ different assets, but that the covariance matrix of the asset returns \emph{does not depend} on the auxiliary features. Therefore, we consider a model  with a constraint that bounds the overall variance of the portfolio. Specifically, if $\Sigma \in \bbR^{d \times d}$ denotes the (positive semidefinite) covariance matrix of the asset returns and $\gamma \geq 0$ is the desired bound on the overall variance (risk level) of the portfolio, then the feasible region $S$ in \eqref{poi} is given by $S := \{w : w^T\Sigma w \leq \gamma, e^Tw \leq 1, w \geq 0\}$. Here $e$ denotes the vector of all ones and since we only require that $e^Tw \leq 1$, the cost vector $c$ in \eqref{poi} represents the negative of the incremental returns of the assets above the risk-free rate. In other words, it holds that $c = -\tilde r$ where $\tilde r = r - r_{\text{RF}}e$, $r$ represents the vector of asset returns, and $r_{RF}$ is the risk-free rate. We use the SGD approach (Algorithm \ref{sgd} of Appendix \ref{sec:sgd}) for training the SPO+ model of method (1.). \edit{Training the SPO+ model takes 3 to 5 minutes for each ERM instance, while the other methods typically take less than a second.}  For brevity, we defer the details of the experimental setup to Appendix \ref{portfolio_experiment}.

\begin{figure}[h!]
\centering
\includegraphics[width=\linewidth]{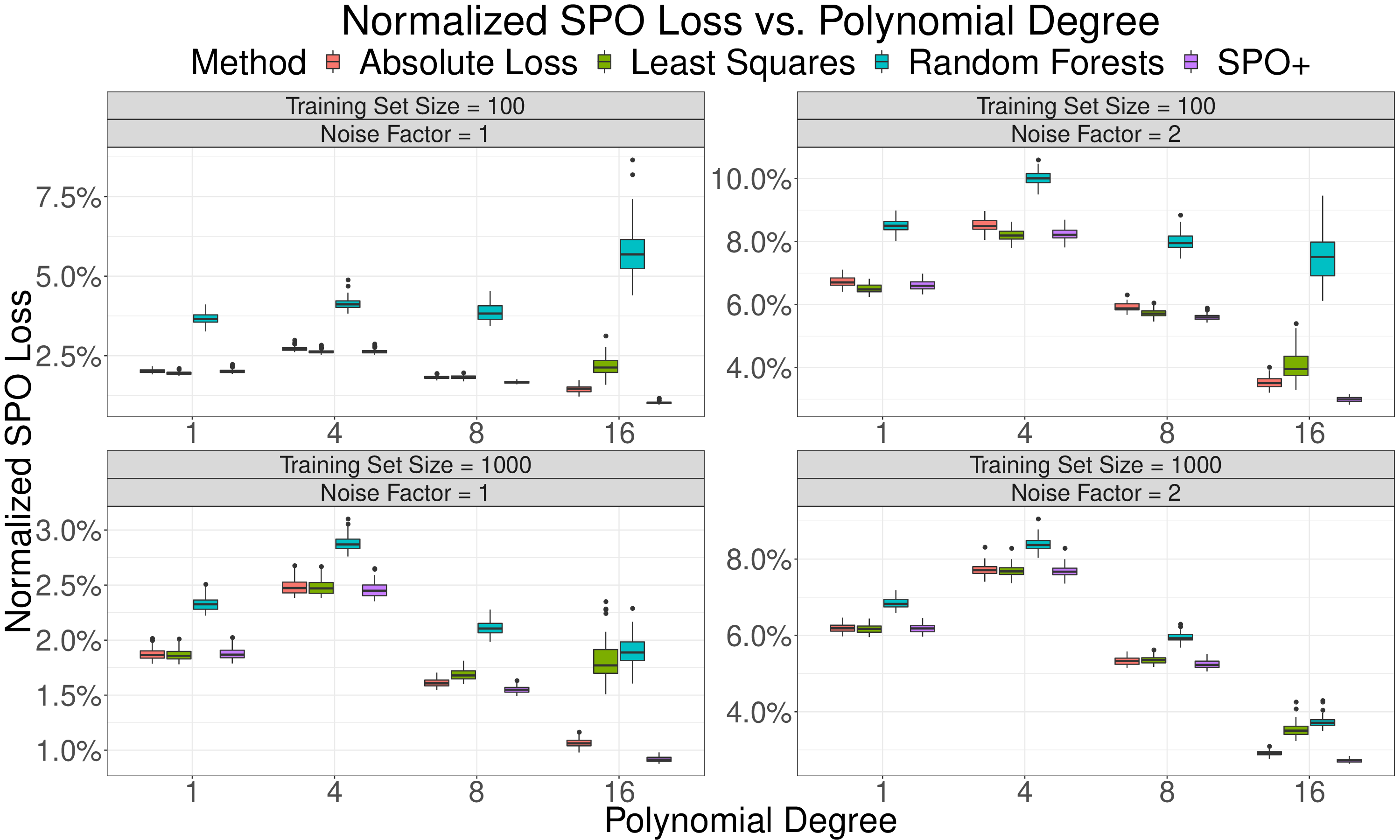}
\caption{Normalized test set SPO loss for the SPO+, least squares, absolute loss, and random forests methods on portfolio optimization instances.}
\label{fig:portfolio_expts}
\end{figure}

Figure \ref{fig:portfolio_expts} displays our results for this experiment. Generally we observe similar patterns as in the shortest path experiment, although comparatively larger values of $\degg$ are needed to demonstrate the relative superiority of SPO+. In summary across all of our experiments, our results indicate that as long as there is some degree of model misspecification, then SPO+ tends to offer significant value over competing approaches, and this value is further strengthened in cases where more data available. The SPO+ approach is either always close to the best approach, or dominating all other approaches, making it a fairly suitable choice across all parameter regimes.

\section{Conclusion}
In this paper, we provide a new framework for developing prediction models under the predict-then-optimize paradigm. Our SPO framework relies on new types of loss functions that explicitly incorporate the problem structure of the optimization problem of interest. Our framework applies for any problem with a linear objective, even when there are integer constraints. 

Since the SPO loss function is nonconvex, we also derived the convex SPO+ loss function using several logical steps based on duality theory. Moreover, we prove that the SPO+ loss is consistent with respect to the SPO loss, which is a fundamental property of any loss function. In fact, our results also directly imply that the least squares loss function is also consistent with respect to the SPO loss. Thus, least squares performs well when the ground truth is near linear, although, at least empirically, SPO+ strongly outperforms all approaches when there is model misspecification. In subsequent work, we have shown how to train decision trees with SPO loss \citep{elmachtoub2020decision} and developed generalization bounds of the SPO loss function \citep{el2019generalization}. Naturally, there are many important directions to consider for future work including more empirical testing and case studies, handling unknown parameters in the constraints, and dealing with nonlinear objectives.

\section*{Acknowledgements}
 The authors gratefully acknowledge the support of NSF Awards CMMI-1763000, CCF-1755705, and CMMI-1762744.

\bibliographystyle{ormsv080} 
\bibliography{SmartPredictions_references} 
\newpage
\begin{APPENDIX}{}
\normalsize
\section{Examples} \label{sec:examples}
\begin{example}[Network Flow] \label{nf}  An example of the nominal optimization problem is a minimum cost network flow problem, where the decisions are how much flow to send on each edge of the network. We assume that the underlying graph is provided to us, e.g., the road network of a city. The feasible region $S$ represents flow conservation, capacity, and required flow constraints on the underlying network.
The cost vector $c$ is not known with certainty, but can be estimated from data $x$ which can include features of time, day, edge lengths, most recent observed cost, and so on. An example of the hypothesis class is the set of linear prediction models given by $\calH = \{f : f(x) = Bx \text{ for some } B \in \bbR^{d \times p}\}$. The linear model can be trained, for example, according to the mean squared error loss function, i.e., $\ell(\hat c, c) = \tfrac{1}{2}\|\hat c - c\|_2^2$. The corresponding empirical risk minimization problem to find the best linear model $B^*$ then becomes $
\min_{B} \ \frac{1}{n}\sum_{i = 1}^n \frac{1}{2}\|B x_i - c_i \|_2^2 $.
The decision rule to find the optimal network flow given a feature vector $x$ is  $w^*(B^*x)$. 
\Halmos
\end{example}

\begin{example}[Necessity of Assumption \ref{assumption-fun}.1] \label{normal-example}
Suppose that $d = 1$, $S = [-1/2, +1/2]$, and $c$ is normally distributed with mean $\bar c = 0$ and variance 1. Then $W^\ast(\bar c) = S$ and $\ellT(\bar c, c) = \supp(c) - z^\ast(c) = |c|$ for all $c$. Clearly though, $\ellT(1, c) = -c/2 + |c|/2$. Moreover, $\ellT(1, c)$ strictly dominates $\ellT(\bar c, c)$ in the sense that $\ellT(1, c) \leq \ellT(\bar c, c)$ for all $c$ and $\ellT(1, c) < \ellT(\bar c, c)$ for $c > 0$. Therefore, $\riskT(1) < \riskT(\bar c)$ and hence $\bar c$ is not a minimizer of $\riskT(\cdot)$. \Halmos
\end{example}

\begin{example}[Necessity of Assumption \ref{assumption-fun}.2]\label{bad_example2}
Define $S$ as the two-dimensional unit square with extreme points at $(0,0), (1,0)$, $(1,1)$ and $(0,1)$. Let $c$ have a continuous distribution on $\mathbb{R}^2$ defined in the following manner. The probability that $c$ occurs in each quadrant is exactly $0.25$. The distribution over each quadrant is a folded normal. The mean of the folded normals is $(9,9)$ in quadrant 1, $(-1,1)$ in quadrant 2, $(-1,-1)$ in quadrant 3, and $(1,-1)$ in quadrant 4. Thus, $\bar{c}= (2,2)$ and $W^{\ast}(\bar{c})=\{(0,0)\}$. Note that all cost vectors in quadrant 1 are minimized by $(0,0)$, quadrant 2 by $(1,0)$, quadrant 3 by $(1,1)$, and quadrant 4 by $(0,1)$. Therefore, $\bbE_c [w^{\ast}(c)]=(0.5,0.5)$. 

Now we claim $c^{\ast}=(0,0)$ is the unique minimizer of SPO+ risk. Since  $w^\ast(c)$ only depends on which quadrant $c$ lies in and both $c$ and $-c$ lie in each of the four quadrants with equal probability, it holds that $\bbE_c[w^{\ast}(2c^{\ast}-c)]=\bbE_c[w^{\ast}(-c)]=\bbE_c[w^{\ast}(c)]=(0.5, 0.5)$ which implies that $c^{\ast}$ is optimal.  Finally, observe that $W^{\ast}(c^{\ast}) = S \not\subseteq \{(0,0)\} =  W^{\ast}(\bar c)$, which means that by Proposition \ref{true_consistency}, $c^{\ast}$ is not a minimizer of the SPO risk.\Halmos
\end{example} 

\begin{example}[Necessity of Assumption \ref{assumption-fun}.3]\label{bad_example1}
Define $S$ as the two-dimensional simplex with extreme points at $(0,0), (1,0),$ and $(0,1)$. Let $c$ have a two point distribution on the points $(-2,1)$ and $(0,-1)$, each with probability of 0.5. One can confirm that $\bar{c}=(-1,0)$, $W^{\ast}(\bar{c})=\{ (1,0) \}$, $\bbE_c[w^{\ast}(c)]=(0.5,0.5)$, and that $c$ is symmetric around its mean ($c$ is equal in distribution to $2 \bar{c}-c$). 

Now we claim $c^{\ast}=(-0.25,-0.5)$ is a minimizer of SPO+ risk. This is confirmed by checking that $\bbE_c[w^{\ast}(2c^{\ast}-c)]=\bbE_c[w^{\ast}(c)]=(0.5, 0.5)$. However, $W^{\ast}(c^{\ast}) = \{(0,1)\} \not\subseteq \{(1,0)\} = W^{\ast}(\bar{c})$. Thus by Proposition \ref{true_consistency}, $c^{\ast}$ is not a minimizer of the SPO risk. \Halmos 
\end{example}

\begin{example}[Necessity of Assumption \ref{assumption-fun}.4] \label{interior}
Let $S \subseteq \bbR^2$ be the interval $S := \{(w_1, w_2)^T : -1 \leq w_1 \leq 1, w_2 = 0\}$ and let $c$ be distributed according to a standard 2-dimensional Gaussian random vector with mean $(0, 0)^T$ (and covariance matrix equal to the identity).
Then, one can easily show that $\riskS(\hat{c}) = \bbE_c[|c_1 - 2\hat{c}_1|] + \sqrt{\frac{2}{\pi}}$. 
Thus, it is apparent that $(0,\beta)^T$ is a minimizer for SPO+ risk for any $\beta \in \mathbb{R}$. Intuitively, this follows from the fact that the second dimension is irrelevant for the optimization problem. \Halmos
\end{example}

\section{Omitted Proofs} \label{sec:proofs}

\subsection{Proof of Proposition \ref{binary1}}
\proof{Proof.}
Let $d = 1$ and the feasible region be the interval $S = [-1/2, +1/2]$. Here the ``cost vector'' $c$ corresponds to a binary class label, i.e., $c$ can take one of two possible values, $-1$ or $+1$. (However, the predicted cost vector is allowed to be any real number.) Notice that, for both possible values of $c$, it holds that $z^\ast(c) = -1/2$. There are three cases to consider for the prediction $\hat c$:  {\em (i)} if $\hat c > 0$ then $W^\ast(\hat c) = \{-1/2\}$ and $\ellT(\hat c, c) = (1 - c)/2$, {\em (ii)} if $\hat c < 0$ then $W^\ast(c) = \{+1/2\}$ and $\ellT(\hat c, c) = (1 + c)/2$, and {\em (iii)} if $\hat c = 0$ then $W^\ast(\hat c) = S = [-1/2, +1/2]$ and $\ellT(\hat c, c) = (1 + |c|)/2 = 1$. Thus, we have $\ellT(\hat c, c) = 0$ when $\hat c$ and $c$ share the same sign, and $\ellT(\hat c, c) = 1$ otherwise. Therefore, $\ellT$ is exactly the 0-1 loss function. \Halmos \endproof

\subsection{Proof of Proposition \ref{alternative-rep}.}
\proof{Proof.}
Let us first prove that the function $q(\alpha) := \supp(c - \alpha \hat c) + \alpha z^\ast(\hat c)$ is monotone decreasing on $\bbR$. Clearly $q(\cdot)$ is a convex function and moreover a subgradient $g$ of $q(\cdot)$ for any $\alpha$ is given by $g := z^\ast(\hat c) - \hat{c}^T w^\ast(\alpha \hat c - c)$. Since $z^\ast(\hat c) = \min\limits_{w \in S} \hat c^T w$, we have that $g \leq 0$ for any $\alpha$. Now, for any $\alpha^\prime \leq \alpha$, the subgradient inequality implies that:
\begin{equation*}
q(\alpha^\prime) ~\geq~ q(\alpha) + g\cdot(\alpha^\prime - \alpha) ~\geq~ q(\alpha) \ ,
\end{equation*}
since $g$ and $\alpha^\prime - \alpha$ are both nonpositive. Thus, $q(\cdot)$ is monotone decreasing. 

Let us now prove that
\begin{equation}\label{stronger-result}
\ellT(\hat c, c) ~=~ \inf_{\alpha \geq 0}\left\{\supp(c - \alpha \hat c) + \alpha z^\ast(\hat c)\right\} - z^\ast(c) \ .
\end{equation}
The proof of \eqref{stronger-result} employs Lagrangian duality (see, e.g., \cite{bertsekas1999nonlinear} and the references therein).
First, note that the set of optimal solutions with respect to $\hat c$, $W^\ast(\hat c) := \arg\min_{w \in S}\left\{\hat c^Tw\right\}$ may be expressed as $W^\ast(\hat c) = S \cap \{w \in \bbR^d : \hat{c}^Tw \leq z^\ast(\hat c)\}$. Therefore, it holds that:
\begin{equation}\label{big_equal}
\begin{array}{clccrl}
\max\limits_{w} & c^Tw & \ \ = \ \ & & \max\limits_{w} & c^Tw \\ 
\text{s.t.} & w \in W^\ast(\hat c) & & & \text{s.t.} &  w \in S \\
& & & & & \hat c^Tw \leq z^\ast(\hat c) \ .
\end{array}
\end{equation}
Let us introduce a scalar Lagrange multiplier $\alpha \in \bbR_+$ associated with the inequality constraint ``$\hat c^Tw \leq z^\ast(\hat c)$'' on the right side of \eqref{big_equal} and then form the Lagrangian:
\begin{equation*}
L(w, \alpha) := c^Tw + \alpha(z^\ast(\hat c) - \hat{c}^Tw) \ .
\end{equation*}
Thus $q(\cdot)$ is exactly the dual function, i.e., it satisfies:
\begin{align*}
\max_{w \in S} \ L(w, \alpha) ~&=~ \max_{w \in S}\left\{c^Tw + \alpha(z^\ast(\hat c) - \hat{c}^Tw)\right\} \\
&=~ \max_{w \in S}\left\{(c - \alpha \hat c)^Tw\right\} + \alpha z^\ast(\hat c) \\
&=~ \xi_S(c - \alpha \hat c) + \alpha z^\ast(\hat c) = q(\alpha) \ .
\end{align*}
Weak duality then implies that $\max\limits_{w \in W^\ast(\hat c)}\left\{c^Tw\right\} ~\leq~ \inf\limits_{\alpha \geq 0}q(\alpha)$ and hence:
\begin{equation*}
\ellT(\hat c, c) ~=~ \max\limits_{w \in W^\ast(\hat c)}\left\{c^Tw\right\} - z^\ast(c) ~\leq~ \inf_{\alpha \geq 0}\left\{\supp(c - \alpha \hat c) + \alpha z^\ast(\hat c)\right\} - z^\ast(c) \ .
\end{equation*}
To prove \eqref{stronger-result}, we demonstrate that strong duality holds by applying Theorem 4.3.8 of \cite{borwein2010convex}. In our setting, the primal problem is the problem on the right-hand side of \eqref{big_equal}. This problem corresponds to the primal minimization problem in \cite{borwein2010convex} by considering the objective function given by $-c^Tw + I_S(w)$ and the constraint function $\hat c^Tw - z^\ast(\hat c)$. (Note that $I_S(w)$ is the convex indicator function equal to $0$ when $w \in S$ and $+\infty$ otherwise.) Since $S$ is a compact and convex set, we satisfy all of the assumptions of Theorem  4.3.8 of \cite{borwein2010convex} and hence strong duality holds. Finally, since $q(\cdot)$ is a monotone decreasing function, we can apply the monotone convergence theorem to replace $\inf_{\alpha \geq 0}$ in \eqref{stronger-result} by $\lim_{\alpha \rightarrow \infty}$.\Halmos  \endproof

{\blockedit
\subsection{Proof of Proposition \ref{prop:convex_ub}.}
\proof{Proof.}
Part (1.) is exactly the chain of inequalities \eqref{der1}-\eqref{der4}. Since the support function $\supp(\cdot)$ is a maximum of linear functions, it is convex and therefore part (2.) follows. To see that part (3.) holds, for any $\tilde{c}$, we have that:
\begin{align*}
\ellS(\tilde c, c) &~=~ \supp(c - 2\tilde c) + 2\tilde{c}^Tw^\ast(c) - z^\ast(c) \\
&~\geq~ \supp(c - 2\hat{c}) + w^\ast(2\hat{c} - c)^T((c - 2\tilde c) - (c - 2\hat c)) + 2\tilde{c}^Tw^\ast(c) - z^\ast(c) \\
&~=~ \supp(c - 2\hat{c}) + 2w^\ast(2\hat{c} - c)^T(\hat{c} - \tilde{c}) + 2\tilde{c}^Tw^\ast(c) - z^\ast(c) \\
&~=~ \supp(c - 2\hat{c}) + 2\hat{c}^Tw^\ast(c) - z^\ast(c) + 2(w^\ast(c) - w^\ast(2\hat{c} - c))^T(\tilde{c} - \hat{c}) \ , 
\end{align*}
where the inequality follows since $w^\ast(2\hat{c} - c) \in \arg\max_{w \in S}\{(c - 2\hat{c})^Tw\} = \partial \supp(c - 2\hat{c})$ (this is a standard result in convex optimization, see, e.g., \cite{boyd2004convex}). Thus, we conclude that $2(w^\ast(c) - w^\ast(2\hat c - c)) \in \partial\ellS(\hat c, c)$. \Halmos \endproof
}

\subsection{Proof of Proposition \ref{binary2}.}

\proof{Proof.}
In the same setup as Proposition \ref{binary1}, we have that $S = [-1/2, +1/2]$ and $c \in \{-1, +1\}$ corresponds to the true label. Note that $\supp(c - 2\hat c) = \tfrac{1}{2}|c - 2\hat c|$, and for $c \in \{-1, +1\}$ we have $w^\ast(c) = -c/2$ and $z^\ast(c) = -1/2$. Therefore, 
$$\ellS(\hat c, c) = \tfrac{1}{2}|c - 2\hat c| + 2\hat{c}(-c/2) + 1/2 = \tfrac{1}{2}|1 - 2\hat{c}c| - \hat{c}c + 1/2 = \max\{0, 1 - 2\hat{c}c\} \ , $$ 
where the second equality follows since $c \in \{-1, +1\}$. \Halmos \endproof

\subsection{Proof of Theorem \ref{thm:main}.}
\proof{Proof.} Let $x \in \calX$ be fixed. Applying Assumptions \ref{assumption-fun}.2, \ref{assumption-fun}.3, and \ref{assumption-fun}.4 with Proposition \ref{spo_consistent} implies that the vector $\E[c|x] = \arg\min\limits_{f(x) \in \bbR^d}\bbE_c\left[\ellS(f(x), c) ~|~ x\right]$, i.e., $\E[c|x]$ is the unique minimizer of $\bbE_c\left[\ellS(f(x), c) ~|~ x\right]$. Thus, since $\mathbb{E}_{(x,c) \sim \calD} [\ellS(f(x), c)] = \bbE_x\left[\bbE_c\left[\ellS(f(x), c) ~|~ x\right]\right]$, any minimizer of the SPO+ risk \eqref{truerisk_spo} must be almost surely equal to the function $x \mapsto \E[c|x]$.

Again, let $x \in \calX$ be fixed. If $W^\ast(\bbE[c | x])$ is a singleton, then Proposition \ref{true_consistency} implies that $\E[c|x] \in \arg\min\limits_{f(x) \in \bbR^d}\bbE_c\left[\ellT(f(x), c) ~|~ x\right]$, i.e., $\E[c|x]$ is a minimizer of $\bbE_c\left[\ellT(f(x), c) ~|~ x\right]$. Since $\mathbb{E}_{(x,c) \sim \calD} [\ellT(f(x), c)] = \bbE_x\left[\bbE_c\left[\ellT(f(x), c) ~|~ x\right]\right]$ and Assumption \ref{assumption-fun}.1 states that $W^\ast(\bbE[c | x])$ is a singleton with probability 1 over $x$, we have that any function that almost surely equals the function $x \mapsto \bbE[c | x]$ is a minimizer of the SPO risk \eqref{truerisk_true}. In particular, as argued in the preceding paragraph, any minimizer of the SPO+ risk \eqref{truerisk_spo} must be almost surely equal to the function $x \mapsto \E[c|x]$ and therefore minimizes the SPO risk. Thus, we have shown that the SPO+ loss function is Fisher consistent with respect to the SPO loss. \Halmos \endproof

\subsection{Proof of Proposition \ref{true_consistency}.}
\proof{Proof.}
Consider a cost vector $c^\ast$ that is a minimizer of $\riskT(\cdot)$.	Let $\bar w$ be an optimal solution of $P(\bar c)$, i.e., $\bar w \in W^\ast(\bar c)$, and let $\tilde c$ be chosen such that $\bar w$ is the unique optimal solution of $P(\tilde c)$, i.e., $W^\ast(\tilde c) = \{\bar w\}$. (Note that if $\bar w$ is the unique optimal solution of $P(\bar c)$ then it suffices to select $\tilde c \gets \bar c$, otherwise we may take $\tilde c$ as a slight perturbation of $\bar c$). Then it holds that:
\begin{align*}
0 &\leq \riskT(\tilde{c}) - \riskT(c^\ast) ~=~ \bbE_c\left[\max_{w \in W^\ast(\tilde c)}\left\{c^Tw\right\}\right] - \bbE_c \left[\max_{w \in W^\ast(c^\ast)}\left\{c^Tw\right\}\right] \\
&= \bbE_c [c^T\bar w] - \bbE_c \left[\max_{w \in W^\ast(c^\ast)}\left\{c^Tw\right\}\right]
~=~ \bar{c}^T\bar{w} - \bbE_c \left[\max_{w \in W^\ast(c^\ast)}\left\{c^Tw\right\}\right] ~\leq~ \bar{c}^T\bar{w} - \max_{w \in W^\ast(c^\ast)}\left\{\bar{c}^Tw\right\} \ ,
\end{align*}
where the first inequality and first equality are by definition, the second equality uses $W^\ast(\tilde c) = \{\bar w\}$, the third equality uses linearity of expectation, and the final inequality is Jensen's inequality.
Finally, we conclude that, for any $w \in W^\ast(c^\ast)$, it holds that $\bar{c}^Tw \leq \bar{c}^T\bar{w} = z^\ast(\bar{c})$. Therefore, $W^\ast(c^\ast) \subseteq W^\ast(\bar c)$. 

To prove the other direction, consider a cost vector $c^\ast$ such that $W^\ast(c^\ast) = \{w^\ast(c^\ast)\}$ is a singleton and $W^\ast(c^\ast) \subseteq W^\ast(\bar c)$, i.e., $w^\ast(c^\ast) \in W^\ast(\bar c)$. Let $c^{\ast \ast}$ be an arbitrary minimizer of $\riskT(\cdot)$. Then,
\begin{align*}
\riskT(c^\ast)  -  \riskT(c^{\ast \ast}) &= \bbE_c \left[\max_{w \in W^\ast(c^{\ast})}\left\{c^Tw\right\}\right] - \bbE_c \left[\max_{w \in W^\ast(c^{\ast \ast})}\left\{c^Tw\right\}\right] \\
&= \bar{c}^Tw^\ast(c^\ast) - \bbE_c \left[\max_{w \in W^\ast(c^{\ast \ast})}\left\{c^Tw\right\}\right] \\ &= z^\ast(\bar c) - \bbE_c \left[\max_{w \in W^\ast(c^{\ast \ast})}\left\{c^Tw\right\}\right] \\
&\leq z^\ast(\bar c) - \max_{w \in W^\ast(c^{\ast \ast})}\left\{\bar{c}^Tw\right\} ~\leq~ 0 \ ,
\end{align*}
where the second equality uses $W^*(c^\ast) = \{w^\ast(c^\ast)\}$, the third equality uses $w^\ast(c^\ast) \in W^\ast(\bar c)$, and the first inequality is Jensen's inequality.
Finally, we conclude that since $c^{\ast \ast}$ is a minimizer of $\riskT(\cdot)$ and $c^{\ast}$ has at most the same risk, then $c^{\ast}$ is also a minimizer of $\riskT(\cdot)$. \Halmos \endproof

\subsection{Proof of Proposition \ref{spo_consistent}.}
\proof{Proof.}
	
\textit{(a)} First, note that it is straightforward to show that $\riskS(\hat c)$ is finite valued for all $\hat c$, which follows since $S$ is compact and $\bar c := \bbE[c]$ is finite.
Moreover, by Proposition \ref{prop:convex_ub}, it is clear that $\riskS(\cdot)$ is convex on $\bbR^d$. In particular, for any point $\hat c$ the subdifferential $\partial \riskS(\hat c)$ is nonempty and, since $\riskS(\hat c)$ is finite, we have that $\partial \riskS(\hat c) = \bbE_c\left[\partial \ellS(\hat c, c)\right]$ (see \cite{strassen1965existence, bertsekas1973stochastic}), where the latter refers to the selection expectation over the random set $\partial \ellS(\hat c, c)$ (see Section 2.1.2 in \cite{molchanov2005theory}).
By the linearity of selection expectation, note that $\partial \riskS(\hat c) = -2\bbE_c\left[W^\ast(2\hat c - c)\right] + 2\bbE_c[w^\ast(c)] = -2\bbE_c\left[w^\ast(2\hat c - c)\right] + 2\bbE_c[w^\ast(c)]$, where the second equality follows since the distribution of $c$ is continuous on all of $\bbR^d$ which implies that $W^\ast(2\hat c - c)$ is a singleton with probability 1 (see, e.g., the introductory discussion in \cite{drusvyatskiy2011generic}).

Now, the optimality conditions for the convex problem $\min\limits_{\hat c \in \bbR^d} \riskS(\hat c)$ state that $c^\ast$ is a global minimizer if and only if $0 \in \partial \riskS(c^\ast)$. 
By the discussion in the previous paragraph, the optimality conditions may be equivalently written as $\bbE_c [w^\ast(c)] = \bbE_c[w^\ast(2c^\ast - c)]$. 
Finally, since $c$ is centrally symmetric around its mean, we have that $c$ is equal in distribution to $2\bar{c}-c$; hence $\bbE_c [w^\ast(c)]= \bbE_c [w^\ast(2\bar{c}-c)]$. Therefore $\bar c$ satisfies $0 \in \partial \riskS(\bar{c})$ and is an optimal solution of $\min\limits_{\hat c \in \bbR^d} \riskS(\hat c)$.
	
\textit{(b)} Consider any vector $\Delta \neq 0$, and let us rewrite the difference $\riskS(\bar{c} + \Delta) - \riskS(\bar{c})$ as follows:
\begin{equation*}
\begin{array}{ll}
\riskS(\bar{c} + \Delta) -~ \riskS(\bar{c}) & =~  \bbE_c[\supp(c - 2(\bar c + \Delta)) ~+~ 2(\bar c + \Delta)^Tw^\ast(c) ~-~ z^\ast(c)] ~- \\
& \ \ \ \ \bbE_c[\supp(c - 2\bar c) ~+~ 2\bar c^Tw^\ast(c) ~-~ z^\ast(c)] \\
& =~ \bbE_c[\supp(c - 2\bar c - 2\Delta) ~-~ \supp(c - 2\bar c) ~+~ 2\Delta^Tw^\ast(c)] \\
& =~ \bbE_c[\supp(c - 2\bar c - 2\Delta) ~-~ (c - 2\bar c)^Tw^\ast(2\bar c - c) ~+~ 2\Delta^Tw^\ast(c)] \\
& =~ \bbE_c[\supp(c - 2\bar c - 2\Delta) ~-~ (c - 2\bar c - 2\Delta)^Tw^\ast(2\bar c - c)] \\
& =~ \bbE_c[(c - 2\bar c - 2\Delta)^T(w^\ast(2\bar{c} + 2\Delta - c) ~-~ w^\ast(2\bar{c} - c))] \ . \\
\end{array}
\end{equation*}
The first equality follows from the definitions of $\riskS(\cdot)$ and $\ellS(\cdot, \cdot)$, the second follows from linearity of expectation, and the third follows from the definition of $\supp(\cdot)$. The fourth equality follows from the fact that $c$ is symmetric about its mean so $\bbE_c [w^\ast(c)]= \bbE_c [w^\ast(2\bar{c}-c)]$ and again uses linearity of expectation. The fifth equality follows from the definition of $\supp(\cdot)$. 

Since $w^\ast(2\bar{c} + 2\Delta - c)$ is the maximizer for $c - 2\bar c - 2\Delta$, we have that $(c - 2\bar c - 2\Delta)^T(w^\ast(2\bar{c} + 2\Delta - c) - w^\ast(2\bar{c} - c))$ is a nonnegative random variable. Moreover, since the distribution of $c$ is continuous on all of $\bbR^d$ and the interior of $S$ is non-empty, there must be some probability mass where $w^\ast(2\bar{c}+2\Delta-c)$ disagrees with $w^\ast(2\bar{c}-c)$, i.e., it holds that $\bbP(w^\ast(2\bar{c}+2\Delta-c) \neq w^\ast(2\bar{c}-c)) > 0$. Note also that $w^\ast(2\bar{c}+2\Delta-c)$ is the unique maximizer for $c - 2\bar c - 2\Delta$ with probability one. Thus, we have $\bbP((c - 2\bar c - 2\Delta)^T\left(w^\ast(2\bar{c} + 2\Delta - c) - w^\ast(2\bar{c} - c)) > 0\right) > 0$ and therefore
\begin{equation*}
\riskS(\bar{c} + \Delta) ~-~ \riskS(\bar{c}) = \bbE_c[(c - 2\bar c - 2\Delta)^T(w^\ast(2\bar{c} + 2\Delta - c) ~-~ w^\ast(2\bar{c} - c))] ~> 0 \ .
\end{equation*}
Thus $\bar{c} + \Delta$ is not a minimizer of the SPO+ risk, and $\bar{c}$ must be the unique minimizer.
\Halmos \endproof

\subsection{Proof of Proposition \ref{reformulation}}
\proof{Proof.}
Recall that $\ellS(Bx_i, c_i) = \supp(c_i - 2Bx_i) + 2(Bx_i)^Tw^\ast(c_i) - z^\ast(c_i)$. Linear programming strong duality implies that:
\begin{equation*}
\begin{array}{cllccll}
\supp(c_i - 2Bx_i) \ \ = \ \ & \max\limits_{w} & (c_i - 2Bx_i)^Tw & \ \ = \ \ & \min\limits_{p} & -b^Tp_i \\
& \text{s.t.} & Aw \geq b & & \text{s.t.} & -A^Tp = c_i - 2Bx_i \\
& & & & & p \geq 0 \ .
\end{array}
\end{equation*}
Note also that $2(w^\ast(c_i)x_i^T) \bullet B$ is just a rewriting of $2(Bx_i)^Tw^\ast(c_i)$ as an explicit linear function. Thus, introducing variables $p_i \in \bbR^m$ for each $i \in \{1, \ldots, n\}$, it is clear that \eqref{reform} is equivalent to \eqref{erm_true3}. \Halmos \endproof

\section{Stochastic Gradient Approach} \label{sec:sgd}
{\blockedit
The idea of this approach is to apply a stochastic gradient method (see, for example, \citep{RM51}, \cite{bottou2018optimization}, and the references therein) directly to problem \eqref{erm_true3} in its original format. For this approach, we assume that the regularizer is convex and subdifferentiable (i.e., we can efficiently compute subgradients of $\Omega(\cdot)$). For example, the ridge penalty $\Omega(B) = \tfrac{1}{2}\|B\|_F^2$ and $\ell_1$ penalty $\Omega(B) = \|B\|_1$ satisfy this assumption. While the algorithm presented herein is generic and can be applied to any subdifferentiable regularizer, it is worth mentioning that an alternative approach to differentiating the regularizer $\Omega(\cdot)$ is to apply a proximal gradient type method \cite{parikh2014proximal}, which instead relies on utilizing the specific functional form of certain regularizers.   



Let us now describe some more details about how to apply this approach.
For convenience, let us write the objective function of \eqref{erm_true3} as $\LSPO(B) := \tfrac{1}{n}\sum_{i = 1}^n \phi_i(B)$, where $\phi_i(B) := \ellS(Bx_i, c_i) + \lambda\Omega(B)$. As mentioned in Proposition \ref{prop:convex_ub}, $\ellS(\cdot, c)$ is convex (but not necessarily differentiable) for a fixed $c$ and therefore $\phi_i(\cdot)$ is a convex function. The following proposition describes how to compute subgradients of $\phi_i(\cdot)$ and also builds on Proposition \ref{prop:convex_ub}.

\begin{proposition}[Subgradient for SPO+ Linear Regression] \label{prop:phi_subgrad}
Let $(x_i, c_i)$ for $i \in \{1, \ldots, n\}$, $B \in \bbR^{d \times p}$, and $\Psi \in \partial \Omega(B)$ be given. Then, it holds that $2(w^\ast(c_i) - w^\ast(2Bx_i - c_i))x_i^T + \lambda\Psi \in \partial \phi_i(B)$.
\end{proposition}
\proof{Proof.}
Recall from Proposition \ref{prop:convex_ub} that $2(w^\ast(c_i) - w^\ast(2\hat c - c_i)) \in \partial\ellS(\hat c, c_i)$ for any $\hat c$. Thus, letting $\hat c \gets Bx_i$ yields for any $B^\prime$ that:
\begin{align*}
\ellS(B^\prime x_i, c_i) ~&\geq~ \ellS(Bx_i, c_i) ~+~ 2(w^\ast(c_i) - w^\ast(2Bx_i - c_i))^T(B^\prime x_i - Bx_i) \\
&=~ \ellS(Bx_i, c_i) + (2(w^\ast(c_i) - w^\ast(2Bx_i - c_i))x_i^T) \bullet (B^\prime - B) \ ,
\end{align*}
and hence $2(w^\ast(c_i) - w^\ast(2Bx_i - c_i))x_i^T \in \partial \ellS(Bx_i, c_i)$. Therefore, $2(w^\ast(c_i) - w^\ast(2Bx_i - c_i))x_i^T + \lambda\Psi \in \partial\ellS(Bx_i, c_i) + \partial\Omega(B) = \partial\phi_i(B)$.
\Halmos 
\endproof
}

Algorithm \ref{sgd} presents the application of stochastic subgradient descent with mini-batching to problem \eqref{erm_true3}. Algorithm \ref{sgd} is a standard application of stochastic subgradient descent in our setting, and closely follows \cite{nemirovski2009robust}.
It is important to emphasize that the main computational requirement of Algorithm \ref{sgd} is access to the optimization oracle $w^\ast(\cdot)$ for problem \eqref{poi}, which is utilized $N$ times during each iteration of the method. 
\edit{Therefore, Algorithm \ref{sgd} may be applied in any case when a practically efficient optimization oracle is available, including some combinatorial and mixed-integer problems as mentioned in Remark \ref{remark:mixed}.}
The main parameters that need to be set are the batch size parameter $N$ and the step-size sequence $\{\gamma_t\}$ (in addition to the regularization parameter $\lambda$). In Section \ref{sec:computation}, we describe precisely how we set these parameters for our experiments. In general, we recommend setting the batch size parameter to a fixed constant such as $5$ or $10$. The choice of the step-size depends on the properties of the regularizer $\Omega(\cdot)$. For general convex $\Omega(\cdot)$, or simply when $\lambda = 0$, we recommend following \cite{nemirovski2009robust} where it is suggested to set the step-size sequence to $\gamma_t = \tfrac{\theta}{\sqrt{t+1}}$ for a fixed constant $\theta > 0$. 
In the case of the ridge penalty $\Omega(B) = \tfrac{1}{2}\|B\|_F^2$ and when $\lambda > 0$, since this function is strongly convex one may alternatively set the step-size sequence to $\gamma_t = \tfrac{2}{\lambda(t+2)}$. Since the SPO+ loss function is convex but non-smooth (and also Lipschitz continuous), these are essentially the only two options that will lead to a precise convergence rate guarantee. Indeed, for the sequence $\gamma_t = \tfrac{\theta}{\sqrt{t+1}}$, \cite{nemirovski2009robust} derive an $O(1/\epsilon^2)$ complexity guarantee in terms of the number of iterations required to ensure that the averaged iterate is $\epsilon$-suboptimal in expectation. It is also straightforward to see that this guarantee can be extended to cases when the optimization oracle $w^\ast(\cdot)$ is only computed approximately with an additive approximation error, as long as we only desire convergence on the order of the accuracy of the oracle.
On the other hand, \cite{lacoste2012simpler} derive a similar $O(1/(\lambda \epsilon))$ complexity guarantee for the $\gamma_t = \tfrac{2}{\lambda(t+2)}$ sequence.

{\blockedit
\begin{algorithm}
	\caption{Stochastic Subgradient Descent with Mini-Batching for Problem \eqref{erm_true3}}\label{sgd}
	\begin{algorithmic}
		\STATE Initialize  $B_0 \in \bbR^{d \times p}$ (typically $B_0 \gets 0$), $t \gets 0$. Set batch size parameter $N \geq 1$. \\
		$ \ $ \\
		At iteration $t \geq 0$:
		\STATE 1. For $j =1,\ldots,N$: \\
        \ \ \ \ \ \ \ \ \ \ Sample $i$ uniformly at random from the set $\{1, \ldots, n\}$ .\\
		\ \ \ \ \ \ \ \ \ \ Compute $\tilde{w}_t^j \gets w^\ast(2B_tx_i - c_i)$ .\\
        \ \ \ \ \ \ \ \ \ \ Set $\tilde G_t^j \gets (w^\ast(c_i) - \tilde{w}_t^j)x_i^T$ .
		\STATE 2. Select $\gamma_t > 0$ and compute: \\
        \ \ \ \ \ \ \ \ \ \ $\Psi_t \in \partial \Omega(B_t)$ \\
        \ \ \ \ \ \ \ \ \ \ $G_t \gets \tfrac{1}{N}\sum_{j = 1}^N\tilde{G}_t^j ~+~ \lambda\Psi_t$ \\
        \ \ \ \ \ \ \ \ \ \ $B_{t+1} \gets B_t - \gamma_t G_t$ \\
        \ \ \ \ \ \ \ \ \ \ $\bar{B}_t \gets \tfrac{1}{\sum_{s = 0}^t\gamma_{s}} \sum_{s = 0}^t\gamma_{s} B_{s}$ \ .
	\end{algorithmic}
\end{algorithm}
}

It is important to note that Algorithm \ref{sgd} is a standard and basic application of stochastic subgradient descent, and that one may consider making several adjustments to the method in order to improve its practical and possibly theoretical performance (see, e.g., \cite{bottou2012stochastic}, and \cite{bottou2018optimization} and the references therein). 
We also mention that one may employ early stopping with this method, whereby we maintain a validation set to monitor (estimates of) the out-of-sample \emph{true SPO loss} as the algorithm runs and then ultimately choose the averaged iterate $\bar{B}_t$ with the smallest true SPO loss $\ellT$ on the validation set as the final model.


\section{Experimental Details of Portfolio Optimization Application}\label{portfolio_experiment}
In this experiment, we set the number of assets $d = 50$ and the return vectors and features are synthetically generated according to a similar process as before. 
As before, we generate a random matrix $B^\ast \in \bbR^{d \times p}$ that encodes the parameters of the true model, whereby each entry of $B^\ast$ is Bernoulli random variable that is equal to 1 with probability $0.5$.
Then, for some noise level parameter $\tau \geq 0$, we generate a factor loading matrix $L \in \bbR^{50 \times 4}$ such that each entry of $L$ is uniformly distributed on $[-0.0025\tau, 0.0025\tau]$ independently of everything else. 
To generate a training/testing pair $(x_i, c_i)$, as before we first generate the feature vector $x_i \in \bbR^p$ from a multivariate Gaussian distribution with i.i.d. standard normal entries, i.e., $x_i \sim N(0, I_p)$. Then, the incremental return vector $\tilde r_i$ is generated according to the following process:
\begin{enumerate}
\item The conditional mean $\bar r_{ij}$ of the $j^{\text{th}}$ asset return is set equal to $\bar r_{ij} := \left(\tfrac{0.05}{\sqrt{p}}(B^\ast x_i)_j + (0.1)^{1/\degg}\right)^{\degg}$, where $\degg$ is a fixed positive integer parameter.
\item The observed return vector $\tilde r_i$ is set to $\tilde r_i := \bar r_i + Lf + 0.01\tau\varepsilon$, where $f \sim N(0, I_4)$ and $\varepsilon \sim N(0, I_{50})$. The cost vector $c_i$ is set to $c_i := -\tilde r_i$.
\end{enumerate}

It follows from step (2.) above that, conditional on the observed features $x_i \in \bbR^p$, the covariance matrix of the returns is $\Sigma = LL^T + (0.01\tau)^2 I_{50}$. We set the risk level parameter $\gamma$, which is part of the definition of the feasible region $S$, to $\gamma := 2.25 \cdot \bar{w}^T\Sigma\bar{w}$ where $\bar{w} := e/10$ is the equal weight portfolio.
Note that the particular functional form of the conditional mean returns in step (1.) above is also a polynomial function of the features $x$  and that the parameters are set so that $\bar r_{ij} \in [0,1]$ with high probability.
In this experiment, the number of features $p$ is again kept fixed at 5. We varied the training set size $n \in \{100, 1000\}$, the parameter $\degg \in \{1, 4, 8, 16\}$, and the noise level parameter $\tau \in \{1, 2\}$. For each combination of these two parameters, we again ran 50 independent trials and used a test set of size 10,000. We do not use any regularization in this experiment, we use the SGD approach (Algorithm \ref{sgd} of Appendix \ref{sec:sgd}) for training the SPO+ model of method (1.), and in the same manner as before we directly solve the optimization problems of methods (2.) and (3.) using Gurobi. Note that, out of fairness to the other methods, we did not use the early stopping technique for SGD.

\end{APPENDIX}

\end{document}